\documentclass[preprint,12pt]{article}
\usepackage{amssymb}
\usepackage{mathrsfs}
\usepackage{amsfonts}
\usepackage{graphicx}
\usepackage{colortbl,dcolumn}
\usepackage{amsmath}
\usepackage{psfrag}
\usepackage{booktabs}
\usepackage{slashbox}
\numberwithin{equation}{section}
\usepackage[top=1in, bottom=1in, left=1in, right=1in]{geometry}

\newtheorem{thm}{Theorem}[section]

\newtheorem{lem}[thm]{Lemma}
\newtheorem{cor}[thm]{Corollary}
\newtheorem{rem}[thm]{Remark}
\newtheorem{ass}[thm]{Assumption}
\begin{document}

\title{The improved split-step backward Euler method for
stochastic differential delay equations\footnotemark[1]}
       \author{
        Xiaojie Wang \footnotemark[2] \quad  Siqing Gan \footnotemark[3] \\
       {\small School of Mathematical Sciences and Computing Technology,
       Central South University,}\\
      {\small Changsha 410075, Hunan,  PR China } }

       \date{}
       \maketitle

       \footnotetext{\footnotemark[1] This work was supported by NSF of China (No.10871207)
        and the Scientific Research Foundation for the Returned Overseas Chinese
        Scholars, State Education Ministry, and Graduate Research Innovation Projects in Hunan Province (NO.CX2010B118).}
        \footnotetext{\footnotemark[2]Corresponding author: x.j.wang7@gmail.com}
       \footnotetext{\footnotemark[3]Email: siqinggan@yahoo.com.cn}
       \begin{abstract}
          {\rm\small A new, improved split-step backward Euler (SSBE) method is introduced and analyzed for stochastic differential delay equations(SDDEs) with
generic variable delay. The method is proved to be convergent in
mean-square sense under conditions (Assumption \ref{OLC}) that the
diffusion coefficient $g(x,y)$ is globally Lipschitz in both $x$ and
$y$, but the drift coefficient $f(x,y)$ satisfies one-sided
Lipschitz condition in $x$ and globally Lipschitz in $y$. Further,
exponential mean-square stability of the proposed method is
investigated for SDDEs that have a negative one-sided Lipschitz
constant. Our results show that the method has the unconditional
stability property in the sense that it can well reproduce stability
of underlying system, without any restrictions on stepsize $h$.
Numerical experiments and comparisons with existing methods for
SDDEs illustrate the computational efficiency of our method. }\\

%\textbf{PACS: 02.60.Lj}

\end{abstract}

%% main text
\section{Introduction}
In this paper we consider the numerical integration of autonomous
stochastic differential delay equations (SDDEs) in the It\^{o}'s
sense
\begin{equation}
dx(t)=f(x(t),x(t-\tau(t)))dt + g(x(t),x(t-\tau(t)))dw(t)
\label{sddes1}
\end{equation}
with initial data $x(t)=\psi(t), t \in [-\tau,0]$. Here $\tau(t)$ is
a delay term satisfying $\tau(t) \geq 0$ and
$-\tau:=\inf\{t-\tau(t):t\geq 0\}$, $f:\mathbb{R}^d\times
\mathbb{R}^d\longrightarrow\mathbb{R}^d,g:\mathbb{R}^d\times
\mathbb{R}^d\longrightarrow\mathbb{R}^{d\times m}$. We assume that
the initial data is independent of the Wiener measure driving the
equations and $w(t)$ is an $m$-dimensional Wiener process defined on
the complete probability space $(\Omega,\mathcal {F},\{\mathcal
{F}_t\}_{t\geq 0},\mathbb{P})$ with a filtration $\{\mathcal
{F}_t\}_{t\geq 0}$ satisfying the usual conditions (that is, it is
increasing and right continuous while $\mathcal {F}_0$ contains all
$\mathbb{P}$-null sets).

For a given constant stepsize $h>0$, we propose a split-step
backward Euler (SSBE) method for SDDEs (\ref{sddes1}) as follows
\begin{subequations}
\begin{equation}
y_n^* = y_n + hf(y_n^*,\tilde{y}_n^*),\label{ssbe1}
\end{equation}
\begin{equation}
y_{n+1}=y_n^* + g(y_n^*,\tilde{y}_n^*)\Delta w_n, \label{ssbe2}
\end{equation}
\end{subequations}
where $\Delta w_n = w(t_{n+1})-w(t_n)$ and for $0\leq \mu <1,1\leq
q_n\in \mathbb{Z}^+$
\begin{equation}
 \tilde{y}_n^* = \left\{\begin{array}{l}
 \psi(t_n-\tau(t_n)),   \quad  t_n-\tau(t_n) <0,\\
\mu y^*_{n-q_n+1} + (1-\mu)y^*_{n-q_n},  \quad 0\leq
t_n-\tau(t_n)\in[t_{n-q_n},t_{n-q_n+1}) .
\end{array}\right. \label{y*n}
\end{equation}
For arbitrary stepsize $h > 0$, $y_n$ denotes the approximation of
$x(t)$ at time $t_n = nh, n=0,1\cdots$. We remark that $\mu$ in
(\ref{y*n}) depends on how memory values are handled on non-grid
points. Generally there are two ways, the first is to use piecewise
constant interpolation, corresponding to $\mu \equiv0$, and the
second to use piecewise linear interpolation. In later development,
we prefer to assume $0\leq \mu <1$ to cover both cases. Also, we mention
that the scheme (\ref{ssbe1})-(\ref{ssbe2}) here is quite different from the SSBE method in \cite{ZGH09}, which will be explained at the end of this section.

In (\ref{ssbe1})-(\ref{ssbe2}), $y^*_n$ serves as an intermediate
stage value, and in order to continue the process, we have to solve
the implicit equation (\ref{ssbe1}) at every step to acquire
$y^*_n$. Existence and uniqueness of solutions to the implicit
equations (\ref{ssbe1}) will be discussed in section 4. Here, we
always assume that numerical solution of (\ref{ssbe1}) exists
uniquely. And one can easily check that $y^*_n,y_n$ is $
\mathcal{F}_{t_n}$-measurable.

The key aim in this work is to propose a new SSBE method for SDDEs
with variable delay and its convergence and stability in mean-square
sense are investigated under a non-globally Lipschitz condition. This
situation has been investigated in
\cite{HMS02,HMS03,HK05,HK06,HK07,JKN09,YH96,WG09a} for stochastic differential
equations (SDEs) without delay. For SDEs with delay, most of
previous work has been based on the more restrictive assumption that
the coefficients $f,g$ satisfies global Lipschitz and linear growth
conditions, see, for example, \cite{BB00,FMW07,Liu,Mao07,ZGH09}. In
\cite{MS03}, the authors showed that the numerical solution produced
by Euler-Maruyama (EM) method will converge to the true solution of
the SDDEs under the local Lipschitz condition. Note that the proof
of the convergence result in this paper is based on techniques used
in \cite{HMS02,MS03}. In \cite{HMS02}, by interpreting the implicit
method SSBE as the EM applied to a modified SDE the authors were able to
get a strong convergence result. This paper, however, provides an
alternative way to get the convergence result for SSBE. That is, by
giving a direct continuous-time extension we accomplished the
convergence proof for SSBE without considering the modified SDDEs.
Also, in deriving moment bounds of numerical solution, due to the
delay term of our SSBE, i.e., $\tilde{y}_n^*$ in (\ref{ssbe1}),
$y_n^*$ cannot be explicitly dominated by $y_n$ as (3.25) in
\cite{HMS02}. Starting with a recurrence of $y_n^*$ given by
substituting (\ref{ssbe2}) into (\ref{ssbe1}), we overcome this
difficulty and obtained the desired moment bounds. Note that a
similar approach is adopted in the stability analysis.
%new tricks are
%developed to overcome difficulties caused by the delay term of SSBE,
%i.e., $\tilde{y}_n^*$ in (\ref{ssbe1}). More precisely,

Of course, the most important contribution of this work is to
propose an improved SSBE method for SDDEs and to verify its
excellent stability property. In \cite{ZGH09}, the authors proposed
a SSBE method for a linear scalar SDDE with constant lag and its
convergence and stability are studied there. It is worth emphasizing
that our proposed method is a modified version of SSBE in
\cite{ZGH09}. The changes are in two aspects: firstly, we drop the
stepsize restriction $h=\frac{\tau}{\kappa}, \kappa \in
\mathbb{Z}^+$ and allow for arbitrary stepsize $h>0$; secondly and
most importantly, the scheme has been modified to a new one. To see this, the two methods are applied  to a linear scalar SDDE in section \ref{linear_MS}. One can observe that the second terms of $f,g$ in the scheme in \cite{ZGH09} is the numerical solution $y_{n-\kappa+1}$ (see (\ref{SSBEZ}) below). While the corresponding terms in our
scheme is the intermediate stage value $y^*_{n-\kappa}$ (see
(\ref{SSBEW}) below). Note that the modifications of the method do
not raise the strong order of the numerical solution, but they indeed
improve the stability of the method greatly. In fact, it is shown below
that our method can well replicate exponential mean-square stability
of nonlinear test problem, including the linear test equation as a
special case, without any restrictions on stepsize $h$. The convergence and stability results of SSBE can be regarded as an extension of those in \cite{HMS02,HMS03} for SDEs without delay to variable delay case.  This unconditional stability property of (\ref{ssbe1})-(\ref{ssbe2}) demonstrates that the proposed method is promising and will definitely be effective in solving systems with stiffness in the drift term, where
stability investigations are particularly important.
%To our best knowledge, this is
%the very first paper to investigate convergence and stability of
%implicit method for SDDEs under one-sided Lipschitz condition.

%We used a very similar argument in

This article is organized as follows. In next section, a general
convergence result (Theorem \ref{ssbemain}) is established. In
section 3, a convergence result is derived under a one-sided
Lipschitz condition (Assumption \ref{OLC}). Section 4 and 5 are
devoted to exponential mean-square stability property of the method.
%in the
%case where the drift has a negative one-sided Lipschitz constant and
%underlying system has a bounded delay. Positive results are showed
%that the method can well reproduce stability of underlying system,
%without any constraint on stepsize $h$.
Numerical experiments are included in section 6.

%It should also be mentioned that quite a different way from
%\cite{HMS02} is developed to derive the convergence results. In
%\cite{HMS02}, the authors first establish a general convergence
%result for Euler-Maruyama method. Then the authors make efforts to
%adapt the results for EM method to SSBE method and backward Euler
%(BE) method, as SSBE method can be viewed as an Euler-Maruyama
%method applied to a perturbed stochastic differential equations
%(SDEs) of the same form as underlying system. While in our work, we
%give a straight continuous-time extension of the discrete
%approximation and new techniques are developed to accomplish the
%proof of convergence result directly without the help of perturbed
%SDEs.

%\section{Analytical solution of underlying system}

\section{The general convergence results}

Throughout the paper, let $ |\cdot | $ denote both the Euclidean
norm in $ \mathbb{R}^d $ and the trace norm(F-norm) in
$\mathbb{R}^{d\times m}$. As the standing hypotheses, we make the
following assumption.
\begin{ass} \label{LCMC}
The system (\ref{sddes1}) has a unique solution $x(t)$ on
$[-\tau,T]$.  And the functions $f(x,y)$ and $g(x,y)$ are both
locally Lipschitz continuous in $x$ and $y$, i.e., there exists a
constant $L_R$ such that
\begin{equation}
|f(x_2,y_2)-f(x_1,y_1)|^2\vee|g(x_2,y_2)-g(x_1,y_1)|^2\leq
L_R(|x_2-x_1|^2+|y_2-y_1|^2), \label{LLC}
\end{equation} for all $x_1,x_2,y_1,y_2 \in \mathbb{R}^d$
with $|x_1|\vee |x_2|\vee|y_1| \vee |y_2|\leq R$.
\end{ass}

Moreover, we assume that \cite{MS03}
%$\tau = (m-\delta)h, m \geq 1\in \mathcal {N}^+, 0\leq\delta<1.$\\
\begin{ass} \label{init}
$ \psi(t) $ is H\"{o}lder continuous in mean-square with exponent
1/2, that is
\begin{equation}
E|\psi (t) - \psi (s)|^2 \leq \eta_1|t-s|, \label{initial1}
\end{equation}
and $\tau(t)$ is a continuous function satisfying
\begin{equation}
|\tau(t) - \tau(s)| \leq \eta_2|t-s|. \label{initial2}
\end{equation}
\end{ass}
In the following convergence analysis, we find it convenient to use
continuous-time approximation solution. Hence we define continuous
version $\bar{y}(t)$ as follows
\begin{eqnarray}
\bar{y}(t):= \left\{\begin{array}{l}
 \psi(t),   \quad  t\leq0,\\
y_n + (t-t_n)f(y_n^*,\tilde{y}_n^*) + g(y_n^*,\tilde{y}_n^*)\Delta
w_n(t),  \quad t\in [t_n, t_{n+1}), n\geq 0,
\end{array}\right.\label{ce1}
\end{eqnarray}
where $\Delta w_n(t) = w(t)-w(t_n)$. For $t\in [t_n, t_{n+1})$ we
can write it in integral form as follows
\begin{equation}
\bar{y}(t):= y_0 + \int_0^t f(y^*(s),\tilde{y}^*(s))\mbox{d}s +
\int_0^t g(y^*(s), \tilde{y}^*(s))\mbox{d}w_s, \label{ce2}
\end{equation}
where
\begin{equation}
y^*(s) := \sum_{n=0}^{\infty}1_{\{t_n\leq s< t_{n+1}\}} y_n^*, \quad
\tilde{y}^*(s) := \sum_{n=0}^{\infty}1_{\{t_n\leq s< t_{n+1}\}}
\tilde{y}_n^*.
\end{equation}
It is not hard to verify that $\bar{y}(t_n) = y_n$, that is,
$\bar{y}(t)$ coincides with the discrete solutions at the
grid-points.

In additional to the above two assumptions, we will need another
one.
\begin{ass} \label{ANMB}
The exact solution $x(t)$ and its continuous-time approximation
solution $\bar{y}(t)$ have p-th moment bounds, that is, there exist
constants $p>2, A>0$ such that
\begin{equation}
\mathbb{E}\left[\sup\limits_{0 \leq t \leq T}|x(t)|^{p}\right] \vee
\mathbb{E}\left[\sup\limits_{0 \leq t \leq T}|\bar{y}(t)|^{p}\right]
 \leq A \label{MB0}.
\end{equation}
\end{ass}

Now we state our convergence theorem here and give a sequence of
lemmas that lead to a proof.

\begin{thm}\label{ssbemain}Under Assumptions \ref{LCMC},\ref{init},\ref{ANMB}, if
the implicit equation (\ref{ssbe1}) admits a unique solution, then
the continuous-time approximate solution $\bar{y}(t)$ (\ref{ce1})
will converge to the true solution of (\ref{sddes1}) in the
mean-square sense, i.e.,
$$
\mathbb{E} \sup_{0 \leq t \leq T} \left|\bar{y}(t)-x(t)\right|^2
\rightarrow 0, \quad as \quad h \rightarrow 0.
$$
\end{thm}
We need several lemmas to complete the proof of Theorem
\ref{ssbemain}.

First, we will define three stopping times
$$
\rho_R = \inf\{t \geq 0 : |x(t)|\geq R \}, \quad \tau_R = \inf\{t
\geq 0: |\bar{y}(t)|\geq R \:, \mbox{or} \: |y^*(t)|\geq R \}, \quad
\sigma_R = \rho_R \wedge \tau_R,
$$
%$$
%\sigma_R = \rho_R \wedge \tau_R,
%$$
where as usual $\inf \emptyset$ is set as $\infty$ ($\emptyset$
denotes the empty set).

\begin{lem}\label{lem1}
Under Assumption \ref{LCMC}, \ref{init}, there exist
constants
$C_1(R)$, $C_2(R)$ such that for $s \in [t_n, t_{n+1})$ and $h<1$
%\begin{eqnarray}
%&&\mathbb{E}|\bar{y}(s)-y^*(s)|^2 \leq C_1h, \label{yd1}\\
%&&\mathbb{E}1_{\{s\leq\sigma_R\}}|\tilde{y}^*(s)-\bar{y}(s-\tau)|^2
%\leq C_2h. \label{yd2}
%\end{eqnarray}
\begin{eqnarray}
&&\mathbb{E}1_{\{s\leq\sigma_R\}}|\bar{y}(s)-y^*(s)|^2 \leq C_1(R)h, \label{yd1}\\
&&\mathbb{E}1_{\{s\leq\sigma_R\}}|\bar{y}(s-\tau(s))-\tilde{y}^*(s)|^2
\leq C_2(R)h. \label{yd2}
\end{eqnarray}
\end{lem}

{\it Proof.} For $s \in [t_n, t_{n+1})$, by definition of
$\bar{y}(s)$ and $y^*(s)$,
\begin{eqnarray}
\bar{y}(s)-y^*(s) &=& y_n + (s-t_n)f(y_n^*,\tilde{y}_n^*)
+ g(y_n^*,\tilde{y}_n^*)\Delta w_n(s)-y^*_n \nonumber \\
&=& (s-t_{n+1})f(y_n^*,\tilde{y}_n^*) + g(y_n^*,\tilde{y}_n^*)\Delta
w_n(s). \label{dy1}
\end{eqnarray}
Noticing that for $|x|\vee|y|\leq R$
\begin{eqnarray}
|f(x,y)|^2 &\leq& 2|f(x,y)-f(0,0)|^2+2|f(0,0)|^2 \nonumber \\
&\leq& K_R(1+|x|^2+|y|^2), \label{lf}
\end{eqnarray}
with $K_R = 2\max\{L_R,|f(0,0)|\}$. Using linear growth condition of
$g$ and moment bounds in (\ref{MB}), we have appropriate constant
$C_1(R)$ so that
\begin{eqnarray}
\mathbb{E}1_{\{s\leq\sigma_R\}}|\bar{y}(s)-y^*(s)|^2 &\leq&
2K_Rh^2(1+\mathbb{E}|y_n^*|^2+\mathbb{E}|\tilde{y}_n^*|^2)
+2Kh(1+\mathbb{E}|y_n^*|^2+\mathbb{E}|\tilde{y}_n^*|^2) \nonumber \\
&\leq& C_1(R)h. \nonumber
\end{eqnarray}
As for estimate (\ref{yd2}), there are four cases as to the location
of $t_n-\tau(t_n)$ and $s-\tau(s)$:

$\bullet$ 1) $t_n-\tau(t_n) < 0,s-\tau(s) < 0$,

$\bullet$ 2) $t_n-\tau(t_n) \geq 0, s-\tau(s) \geq 0$,

$\bullet$ 3) $t_n-\tau(t_n) < 0, s-\tau(s) \geq 0$,

$\bullet$ 4) $t_n-\tau(t_n) \geq 0, s-\tau(s) < 0$.\\
Noticing that the delay $\tau(s)$ satisfies Lipschitz condition
(\ref{initial2}), one sees that
\begin{equation}
|s-\tau(s) - t_n + \tau(t_n)| \leq (\eta_2 +1)h. \label{tau}
\end{equation}
In the case 1), combining H\"{o}lder continuity of initial data
(\ref{initial1}) and (\ref{tau}) gives the desired assertion. In the
case 2), without loss of generality, we assume $s-\tau(s) \in
[t_i,t_{i+1})$, $t_n-\tau(t_n) = (1-\mu)t_j + \mu t_{j+1} \in
[t_j,t_{j+1}), \: i
> j \geq 0$. Thus we have from (\ref{ssbe1}) and (\ref{y*n2}) that
\begin{eqnarray}
\bar{y}(s-\tau(s))-\tilde{y}^*(s) &=& y_i +
(s-\tau(s)-t_i)f(y^*_i,\tilde{y}^*_i) + g(y^*_i,\tilde{y}^*_i)\Delta
w_i(s-\tau(s)) \nonumber \\ && -(1-\mu) y^*_j -\mu y^*_{j+1} \nonumber  \\
&=& (s-\tau(s)-t_{i+1})f(y^*_i,\tilde{y}^*_i) +
g(y^*_i,\tilde{y}^*_i)\Delta w_i(s-\tau(s))
\nonumber \\ && +(1-\mu)(y_i^*- y^*_j) + \mu(y_i^*-y^*_{j+1}) \allowdisplaybreaks \nonumber \\
&=& (s-\tau(s)-t_{i+1})f(y^*_i,\tilde{y}^*_i) +
g(y^*_i,\tilde{y}^*_i)\Delta w_i(s-\tau(s)) \nonumber \\ && +(1-\mu)
\sum_{k=j}^{i-1} \left[h f(y^*_{k+1},\tilde{y}^*_{k+1})+ g(y^*_{k},
\tilde{y}^*_{k})\Delta w_{k} \right] \nonumber \\ && + \mu
\sum_{k=j+1}^{i-1} \left[h f(y^*_{k+1},\tilde{y}^*_{k+1})+
g(y^*_{k}, \tilde{y}^*_{k})\Delta w_{k} \right], \label{yd5}
\end{eqnarray}
where as usual we define the second summation equals zero when
$i=j+1$. Noticing from (\ref{tau}) that $i-j \leq \eta_2 +1$, and
combining local linear growth bound (\ref{lf}) for $f$, global
linear growth condition for $g$ and moment bounds (\ref{MB}), we can
derive from (\ref{yd5}) that
$$
\mathbb{E}1_{\{s\leq\sigma_R\}}|\tilde{y}^*(s)-\bar{y}(s-\tau(s))|^2
\leq C_2(R)h.
$$
In the case 3) and 4), using an elementary inequality gives
\begin{equation}\nonumber
\begin{split}
&\mathbb{E}1_{\{s\leq\sigma_R\}}|\tilde{y}^*(s)-\bar{y}(s-\tau(s))|^2
\\ \leq& 2\mathbb{E}1_{\{s\leq\sigma_R\}}|\tilde{y}^*(s)-\bar{y}(0)|^2 +
2\mathbb{E}1_{\{s\leq\sigma_R\}}|\bar{y}(0)-\bar{y}(s-\tau(s))|^2.
\end{split}
\end{equation}
Then combining this with results obtained in case 1) and 2) gives
the required result, with $C_2(R)$ a universal constant independent
of $h$.
\begin{lem}\label{lem2}
Under Assumption \ref{LCMC}, \ref{init}, for stepsize $h < 1$, there
exists a constant $C_R$ such that
$$
\mathbb{E} [\sup_{0 \leq t \leq T} |\bar{y}(t \wedge \sigma_R)-x(t
\wedge \sigma_R)|^2 ] \leq C_R h,
$$
with $C_R$ dependent on $R$, but independent of $h$.
\end{lem}

{\it Proof.} For simplicity, denote
$$
e(t) := \bar{y}(t)-x(t).
$$
From (\ref{sddes1}) and (\ref{ce2}), we have
\begin{eqnarray}
&&\mathbb{E} \left[\sup_{0 \leq s \leq t}\left|e(s \wedge
\sigma_R)\right|^2 \right] = \mathbb{E} \left[\sup_{0 \leq s \leq
t}\left|\bar{y}(s\wedge
\sigma_R)-x(s\wedge \sigma_R)\right|^2 \right] \nonumber \\
&=& \mathbb{E} \left[\sup_{0 \leq s \leq t}\left|\int_0^{s\wedge
\sigma_R}f(y^*(r),\tilde{y}^*(r))-f(x(r),x(r-\tau(r)))\mbox{d}r \right.\right.\nonumber \\
&& \left.\left.+ \int_0^{s\wedge
\sigma_R}g(y^*(r),\tilde{y}^*(r))-g(x(r),x(r-\tau(r)))\mbox{d}w(r)\right|^2\right]
\nonumber \\
&\leq&  2T \mathbb{E} \int_0^{t\wedge \sigma_R}
\left|f(y^*(s),\tilde{y}^*(s))-f(x(s),x(s-\tau(s)))\right|^2
\mbox{d}s \allowdisplaybreaks \nonumber
\\ && + 2\mathbb{E}\left[\sup_{0 \leq s \leq t} \left|\int_0^{s\wedge
\sigma_R}g(y^*(r),\tilde{y}^*(r))-g(x(r),x(r-\tau(r)))\mbox{d}w(r)\right|^2\right]
\allowdisplaybreaks \nonumber \\ &\leq& 2(T+4)L_R \mathbb{E}
\int_0^{t\wedge \sigma_R}
|y^*(s)-x(s)|^2+|\tilde{y}^*(s)-x(s-\tau(s))|^2 \mbox{d}s,
\label{3.14}
\end{eqnarray}
where H\"{o}lder's inequality and the Burkholder-Davis-Gundy inequality
were used again. Using the elementary inequality $|a+b|^2 \leq 2|a|^2 +
2|b|^2$, one computes from (\ref{3.14}) that
\begin{eqnarray}
&&\mathbb{E}[\sup_{0 \leq s \leq t}|e(s \wedge \sigma_R)|^2 ]
\nonumber \\ &\leq& 4(T+4)L_R \mathbb{E} \int_0^{t\wedge \sigma_R}
|y^*(s)- \bar{y}(s)|^2 + |\bar{y}(s)-x(s)|^2 \mbox{d}s \nonumber \\
&& + 4(T+4)L_R \mathbb{E} \int_0^{t\wedge \sigma_R} |\tilde{y}^*(s)-
\bar{y}(s-\tau(s))|^2 + |\bar{y}(s-\tau(s))-x(s-\tau(s))|^2 \mbox{d}s \allowdisplaybreaks \nonumber \\
&\leq& 8(T+4)L_R \int_0^t \mathbb{E} [\sup_{0 \leq r \leq
s}|\bar{y}(r\wedge \sigma_R)-x(r\wedge \sigma_R)|^2] \mbox{d}s
\nonumber \\ && + 4(T+4)L_R \mathbb{E} \int_0^{t\wedge \sigma_R}
|y^*(s)-
\bar{y}(s)|^2 \mbox{d}s \nonumber \\
&& + 4(T+4)L_R \mathbb{E} \int_0^{t\wedge \sigma_R} |\tilde{y}^*(s)-
\bar{y}(s-\tau(s))|^2 \mbox{d}s, \label{3.15}
\end{eqnarray}
where the fact was used that $|\bar{y}(s-\tau(s))-x(s-\tau(s))|^2
\leq \sup\limits_{0\leq r\leq s}|\bar{y}(r))-x(r)|^2 $. By taking
Lemma \ref{lem1} into account, we derive from (\ref{3.15}) that,
with suitable constants $\tilde{C}_R, \bar{C}_R$
\begin{eqnarray}
\mathbb{E}[\sup_{0 \leq s \leq t}|e(s \wedge \sigma_R)|^2] &\leq&
8(T+4)L_R \int_0^t \mathbb{E} [\sup_{0 \leq r \leq
s}|\bar{y}(r\wedge \sigma_R)-x(r\wedge \sigma_R)|^2] \mbox{d}s
\nonumber \\
&& + 4(T+4)TL_RC_1(R)h + 4(T+4)TL_RC_2(R)h \nonumber \\
&=& \tilde{C}_R\int_0^t \mathbb{E} [\sup_{0 \leq r \leq s}|e(r
\wedge \sigma_R)|^2] \mbox{d}s + \bar{C}_Rh.
\end{eqnarray}
Hence continuous Gronwall inequality gives the assertion.

{\it Proof of Theorem \ref{ssbemain}.} Armed with Lemma \ref{lem2}
and Assumption \ref{ANMB}, the result may be proved using a similar
approach to that in \cite[Theorem 2.2]{HMS02} and \cite[Theorem
2.1]{MS03}, where under the local Lipschitz condition they showed
the strong convergence of the EM method for the SODEs and SDDEs,
respectively.

\begin{rem}
Under the global Lipschitz condition and linear growth condition (cf
\cite{Mao97}), we can choose uniform constants $C_1(R)$, $C_2(R),
C_R$ in previous Lemma \ref{lem1},\ref{lem2} to be independent of
$R$. Accordingly we can recover the strong order of 1/2 by deriving
$$
\mathbb{E} [\sup_{0 \leq t \leq T} |\bar{y}(t )-x(t)|^2 ] \leq C h,
$$
where $C$ is independent of $R$ and $h$.
\end{rem}

\section{Convergence with a one-sided Lipschitz condition}
In this section, we will give some sufficient conditions on
equations (\ref{sddes1}) to promise a unique global solution of
SDDEs and a well-defined solution of the SSBE method. We make the
following assumptions on the SDDEs.

\begin{ass} \label{OLC}
The functions $f(x,y)$ are continuously differentiable in both $x$
and $y$, and there exist constants $\gamma_1, \gamma_2, \gamma_3,
\gamma_4$, such that $\forall x,y,x_1,x_2,$ $y_1,y_2 \in
\mathbb{R}^d$
\begin{eqnarray}
\langle x_2-x_1, f(x_2,y)-f(x_1,y)\rangle &\leq& \gamma_1|x_2-x_1|^2, \label{OLC1}\\
|f(x,y_2) - f(x,y_1)| &\leq& \gamma_2|y_2-y_1|, \label{OLC2} \\
|g(x_2,y_2) - g(x_1, y_1)|^2 &\leq& \gamma_3|x_2-x_1|^2 +
\gamma_4|y_2-y_1|^2.\label{OLC3}
\end{eqnarray}
\end{ass}
The inequalities (\ref{OLC1}),(\ref{OLC2}) indicate that the first
argument $x$ of $f$ satisfies one-sided Lipschitz condition and the
second satisfies global Lipschitz condition. It is worth noticing
that conditions of the same type as (\ref{OLC1}) and (\ref{OLC2})
have been exploited successfully in the analysis of numerical
methods for deterministic delay differential equations (DDEs)(see
\cite{BZ03} and references therein). As for SDEs without delay, the
conditions (\ref{OLC1}) and (\ref{OLC3}) has been used in
\cite{HMS02,HMS03,HK05,YH96,WG09a}.

We compute from (\ref{OLC1})-(\ref{OLC3}) that
\begin{eqnarray} \label{4.1}
\langle x, f(x,y) \rangle &=& \langle x, f(x,y)-f(0,y)\rangle+
\langle x, f(0,y)-f(0,0) \rangle+\langle x, f(0,0) \rangle \nonumber
\\ &\leq&
(\gamma_1+1)|x|^2+\frac{1}{2}\gamma_2|y|^2+\frac{1}{2}|f(0,0)|^2,
\end{eqnarray}
\begin{equation} \label{4.2}
|g(x,y)|^2 \leq 2|g(x,y)-g(0,0)|^2 + 2|g(0,0)|^2 \leq
2\gamma_3|x|^2+2\gamma_4|y|^2+2|g(0,0)|^2.
\end{equation}
On choosing the constant $K$ as
$$
K=\max \left\{\gamma_1+1,2\gamma_3, \frac{1}{2}\gamma_2, 2\gamma_4,
\frac{1}{2}|f(0,0)|^2, 2|g(0,0)|^2\right\},
$$
the following condition holds
\begin{equation} \label{MC}
x^{T}f(x,y) \vee |g(x,y)|^2 \leq K(1 + |x|^2 + |y|^2), \quad \forall
x,y \in \mathbb{R}^d.
\end{equation}

In what follows we always assume that for $\forall p>0$ the initial
data satisfies
$$
\mathbb{E}\|\psi\|^p := \mathbb{E}\sup_{- \tau \leq s \leq
0}|\psi(s)|^p< \infty.
$$

\begin{thm} \label{EU}
Assume that Assumption \ref{OLC} is fulfilled. Then there exists a
unique global solution $x(t)$ to system (\ref{sddes1}). Morever, for
any $p \geq 2$, there exists constant $C=C(p,T)$
$$
\mathbb{E}\left[\sup\limits_{0 \leq t \leq T}|x(t)|^{p}\right] \leq
C(1+\mathbb{E}\|\psi\|^p).
$$
\end{thm}
{\it Proof.} See the Appendix.

\begin{lem}\label{MBlem}
Assume that $f,g$ satisfy the condition (\ref{MC}) and $h<1$ is
sufficiently small , then for $p\geq 2$ the following moment bounds
hold
\begin{equation}
\mathbb{E}\left[\sup\limits_{0 \leq nh \leq T}|y^*_n|^{2p}\right]
\vee \mathbb{E}\left[\sup\limits_{0 \leq t \leq
T}|\tilde{y}^*(t)|^{2p}\right] \vee \mathbb{E}\left[\sup\limits_{0
\leq t \leq T}|\bar{y}(t)|^{2p}\right]
 \leq A \label{MB}.
\end{equation}
\end{lem}

{\it Proof.} %In (\ref{ssbe1}), replacing $y_n$ by (\ref{ssbe2})
%leads to
%\begin{eqnarray}
%y^*_n = y^*_{n-1} + h f(y^*_n,\tilde{y}^*_n)+
%g(y^*_{n-1},\tilde{y}^*_{n-1})\Delta w_{n-1}, \quad n \geq 1.
%\label{y*n2}
%\end{eqnarray}
%Hence
%$$
%|y^*_n-h
%f(y^*_n,\tilde{y}^*_n)|^2=|y^*_{n-1}+g(y^*_{n-1},\tilde{y}^*_{n-1})\Delta
%w_{n-1}|^2.
%$$
%Expanding it and employing (\ref{MC}) yields
%\begin{eqnarray} \label{3.1}
%|y^*_n|^2 - 2Kh(1+|y^*_n|^2 +|\tilde{y}^*_n|^2) &\leq& |y^*_{n-1}|^2
%+ 2\left \langle y^*_{n-1},g(y^*_{n-1},\tilde{y}^*_{n-1})\Delta
%w_{n-1} \right \rangle \nonumber \\
%&&+|g(y^*_{n-1},\tilde{y}^*_{n-1})\Delta w_{n-1}|^2.
%\end{eqnarray}
Inserting (\ref{ssbe2}) into (\ref{ssbe1}) gives
\begin{eqnarray}
y^*_n = y^*_{n-1} + h f(y^*_n,\tilde{y}^*_n)+
g(y^*_{n-1},\tilde{y}^*_{n-1})\Delta w_{n-1}, \quad n \geq 1.
\label{y*n2}
\end{eqnarray}
Hence
$$
|y^*_n-h
f(y^*_n,\tilde{y}^*_n)|^2=|y^*_{n-1}+g(y^*_{n-1},\tilde{y}^*_{n-1})\Delta
w_{n-1}|^2.
$$
Expanding it and employing (\ref{MC}) yields
\begin{eqnarray} \label{3.1}
|y^*_n|^2 - 2Kh(1+|y^*_n|^2 +|\tilde{y}^*_n|^2) &\leq& |y^*_{n-1}|^2
+ 2\left \langle y^*_{n-1},g(y^*_{n-1},\tilde{y}^*_{n-1})\Delta
w_{n-1} \right \rangle \nonumber \\
&&+|g(y^*_{n-1},\tilde{y}^*_{n-1})\Delta w_{n-1}|^2.
\end{eqnarray}
By definition of $\tilde{y}^*_n$, one obtains $|\tilde{y}^*_n|^2
\leq |y^*_n|^2 + \max_{0 \leq i \leq n-1} |y^*_i|^2 + \|\psi\|^2$.
Taking this inequality into consideration and letting $h<h_0 <
1/(4K)$, we have from (\ref{3.1}) that
\begin{eqnarray}
(1-4Kh)|y^*_n|^2 \leq |y^*_{n-1}|^2 +2Kh(1+\max_{0 \leq i \leq n-1}
|y^*_i|^2 + \|\psi\|^2) \nonumber \\  + 2\left \langle
y^*_{n-1},g(y^*_{n-1},\tilde{y}^*_{n-1})\Delta w_{n-1} \right
\rangle +|g(y^*_{n-1},\tilde{y}^*_{n-1})\Delta w_{n-1}|^2.
\end{eqnarray}
Denoting $\alpha = 1/(1-4Kh_0)$, one computes that
\begin{eqnarray}
|y^*_n|^2 &\leq& |y^*_{n-1}|^2 +6K\alpha h\max_{0 \leq i \leq n-1}
|y^*_i|^2 + 2K\alpha h+ 2K\alpha h \|\psi\|^2 \nonumber \\
&& + 2\alpha\left \langle
y^*_{n-1},g(y^*_{n-1},\tilde{y}^*_{n-1})\Delta w_{n-1} \right
\rangle +\alpha|g(y^*_{n-1},\tilde{y}^*_{n-1})\Delta w_{n-1}|^2.
\end{eqnarray}
By recursive calculation, we obtain
\begin{eqnarray}
|y^*_n|^2 &\leq& |y^*_0|^2 + 6K\alpha h\sum_{j=0}^{n-1}\max_{0 \leq
i \leq j}
|y^*_i|^2 + 2K\alpha T + 2K\alpha T \|\psi\|^2 \nonumber \\
&& + 2\alpha \sum_{j=0}^{n-1} \left \langle
y^*_{j},g(y^*_{j},\tilde{y}^*_{j})\Delta w_{j} \right \rangle
+\alpha \sum_{j=0}^{n-1} |g(y^*_{j},\tilde{y}^*_{j})\Delta w_{j}|^2.
\nonumber
\end{eqnarray}
Raising both sides to the power $p$ gives
\begin{eqnarray}
|y^*_n|^{2p} &\leq& 5^{p-1}\left\{ |y^*_0|^{2p} + (6K\alpha h)^p
n^{p-1}\sum_{j=0}^{n-1}\max_{0 \leq i \leq j}
|y^*_i|^{2p}+\left[2K\alpha T + 2K\alpha T \|\psi\|^2\right]^p \right.\nonumber \\
&& \left.+ (2\alpha)^p \left[\sum_{j=0}^{n-1} \left \langle
y^*_{j},g(y^*_{j},\tilde{y}^*_{j})\Delta w_{j} \right \rangle
\right]^p +\alpha^p n^{p-1} \sum_{j=0}^{n-1}
|g(y^*_{j},\tilde{y}^*_{j})\Delta w_{j}|^{2p} \right \}. \nonumber
\end{eqnarray}
Thus
\begin{eqnarray}
&&\mathbb{E}\max_{1 \leq n \leq M}|y^*_n|^{2p} \nonumber \\ &&\leq
5^{p-1} \left \{\mathbb{E}|y^*_0|^{2p} +(6K\alpha)^p T^{p-1}h
\mathbb{E}\sum_{j=0}^{M-1}\max_{0 \leq i \leq j} |y^*_i|^{2p}\right.
+ \mathbb{E}\left(2K\alpha T + 2K\alpha T \|\psi\|^2\right)^p
\nonumber
\\ && \left.+ (2\alpha)^p\mathbb{E}\max_{1 \leq n \leq M}\left[\sum_{j=0}^{n-1}
\langle y^*_j,g(y^*_j,\tilde{y}^*_j)\Delta w_j\rangle\right]^p +
\alpha^p
M^{p-1}\mathbb{E}\sum_{j=0}^{M-1}|g(y^*_j,\tilde{y}^*_j)\Delta
w_j|^{2p}\right \}.\label{y*n3}
\end{eqnarray}
Here $1\leq M \leq N$, where $N$ is the largest integer number such
that $Nh \leq T$. Now, using the Burkholder-Davis-Gundy inequality
(Theorem 1.7.3 in \cite{Mao97}) gives
\begin{eqnarray}
&&\mathbb{E} \max_{1 \leq n \leq
M}\left[\sum_{j=0}^{n-1}\left\langle
y^*_j,g(y^*_j,\tilde{y}^*_j)\Delta w_j\right \rangle\right]^p \leq
C_p\mathbb{E}\left[
\sum_{j=0}^{M-1}|y^*_j|^2|g(y^*_j,\tilde{y}^*_j)|^2h\right]^{p/2} \allowdisplaybreaks\nonumber\\
&&\leq C_p (Kh)^{p/2}M^{p/2-1}\mathbb{E}\left[\sum_{j=0}^{M-1}
|y^*_j|^p\left(1+|y^*_j|^2+|\tilde{y}^*_j|^2\right)^{p/2}\right] \nonumber\\
&& \leq \frac{1}{2}C_pK^{p/2}T^{p/2-1}h
\mathbb{E}\left[\sum_{j=0}^{M-1} \left(|y^*_j|^{2p} +
3^{p-1}(1+|y^*_j|^{2p}+|\tilde{y}^*_j|^{2p})\right)\right].
\label{estimate1}
\end{eqnarray}
Noticing that
\begin{eqnarray}
\mathbb{E}|\tilde{y}^*_j|^{2p} \leq\mathbb{E} \max_{0\leq i\leq
j}|y_i^*|^{2p}+\mathbb{E}\|\psi\|^{2p} \label{yh},
\end{eqnarray}
inserting it into (\ref{estimate1}), we can find out appropriate
constants $\bar{C}=\bar{C}(p,K,T)$ such that
\begin{eqnarray}
&&\mathbb{E} \max_{0 \leq n \leq
M}\left[\sum_{j=0}^{n-1}\left\langle
y^*_j,g(y^*_j,\tilde{y}^*_j)\Delta w_j\right \rangle\right]^p
\nonumber \\ &\leq& \bar{C}h \sum_{j=0}^{M-1}\mathbb{E}\max_{0 \leq
i \leq j} |y_i^*|^{2p} +\bar{C}(\mathbb{E}\|\psi\|^{2p}+1).
\label{estimate2}
\end{eqnarray}
%where $\check{C} = \frac{1}{2}C_pK^{p/2}T^{p/2-1}(1+2\cdot3^{p-1}),
%\bar{C} = \frac{1}{2}C_pK^{p/2}T^{p/2-1}3^{p-1}(1+\mathbb{E}\|\psi\|^{2p})$.
At the same time, noting the fact $y^*_n,\tilde{y}^*_n\in \mathcal
{F}_{t_n}$ and $\Delta w_n$ is independent of $\mathcal {F}_{t_n}$,
one can compute that, with $\hat{C}=\hat{C}(p,T)$ a constant that
may change line by line
\begin{eqnarray}
\mathbb{E}\sum_{j=0}^{M-1}|g(y^*_j,\tilde{y}^*_j)\Delta w_j|^{2p}
&\leq& \sum_{j=0}^{M-1}\mathbb{E}|g(y^*_j,\tilde{y}^*_j)|^{2p}
\mathbb{E}|\Delta w_j|^{2p} \allowdisplaybreaks
\nonumber \\
&\leq& \hat{C}h^p
\sum_{j=0}^{M-1}\left[1+\mathbb{E}|y^*_j|^{2p}+\mathbb{E}|\tilde{y}^*_j|^{2p}\right]
\allowdisplaybreaks \nonumber \\
&\leq& \hat{C}h^{p-1}(\mathbb{E}\|\psi\|^{2p}+1) + \hat{C}h^p
\sum_{j=0}^{M-1}\mathbb{E}\max_{0 \leq i \leq j} |y^*_i|^{2p}.
\label{estimate3}
\end{eqnarray}
%where $\hat{C} = 3^{p-1}(mK)^pT(1+\mathbb{E}\|\psi\|^{2p}),
%\tilde{C} = 2\cdot 3^{p-1}(mK)^p$.
%Last, we need to bound $\mathbb{E}|y^*_0|^{2p}$. By definition
%(\ref{ssbe1}), one sees that
%$$
%|y^*_0-hf(y^*_0,\tilde{y}^*_0)|^2 = |y_0|^2.
%$$
%Then using a similar approach used before, we can find out a
%constant $\tilde{C} = \tilde{C}(p,T)$ to ensure that
%\begin{equation}
%\mathbb{E}|y^*_0|^{2p} < \tilde{C}. \label{y*0}
%\end{equation}
By definition (\ref{ssbe1}), one sees that
$$
|y^*_0-hf(y^*_0,\tilde{y}^*_0)|^2 = |y_0|^2.
$$
Then using a similar approach used before, we can find out a
constant $c_0 = c_0(p,K)$ to ensure that
\begin{equation}
\mathbb{E}|y^*_0|^{2p} < c_0(\mathbb{E}\|\psi\|^{2p}+1) < \infty.
\label{y*0}
\end{equation}
Inserting (\ref{estimate2}),(\ref{estimate3}) into (\ref{y*n3}) and
considering (\ref{y*0}) and $h<1$ , we have, with suitable constants
$C'=C'(p,K,T),C''=C''(p,K,T)$
\begin{eqnarray} \label{mbend}
\mathbb{E}\max_{0 \leq n \leq M}|y_n^*|^{2p} &\leq&
\mathbb{E}|y^*_0|^{2p} + \mathbb{E}\max_{1 \leq n \leq
M}|y_n^*|^{2p} \nonumber \\ &\leq& C'(\mathbb{E}\|\psi\|^{2p}+1) +
C''h \sum_{j=0}^{M-1}\mathbb{E}\max_{0 \leq i \leq j} |y_i^*|^{2p}.
\end{eqnarray}
Thus using the discrete-type Gronwall inequality, we derive from
(\ref{mbend}) that $\mathbb{E}\left[\sup_{0 \leq nh \leq
T}|y^*_n|^{2p}\right]$ is bounded by a constant independent of $N$.
Then by considering the elementary inequality $|\mu x
+(1-\mu)y|^{2p} \leq \mu|x|^{2p}+(1-\mu)|y|^{2p}$, boundedness of
$\mathbb{E}\left[\sup_{0 \leq t \leq T}|\tilde{y}^*(t)|^{2p}\right]$
is immediate.

To bound $\mathbb{E}\left[\sup_{0 \leq t \leq
T}|\bar{y}(t)|^{2p}\right]$, we shall first bound
$\mathbb{E}\left[\sup_{0 \leq nh \leq T}|y_n|^{2p}\right]$. From
(\ref{ssbe2}), we have
\begin{eqnarray}
\mathbb{E}\left[\sup_{0 \leq nh \leq T}|y_n|^{2p}\right] &\leq&
2^{2p-1}\left\{ \mathbb{E}\left[\sup_{0 \leq nh \leq
T}|y^*_n|^{2p}\right] + \mathbb{E}\left[\sup_{0 \leq nh \leq
T}|g(y^*_n,\tilde{y}_n^*)\Delta w_n|^{2p}\right]\right\} \nonumber
\\ &\leq& 2^{2p-1}\left\{ \mathbb{E}\left[\sup_{0 \leq nh \leq
T}|y^*_n|^{2p}\right] +
\mathbb{E}\sum_{j=0}^{N}|g(y^*_j,\tilde{y}^*_j)\Delta
w_j|^{2p}\right\}. \nonumber
\end{eqnarray}
Now (\ref{estimate3}) and bound of $\mathbb{E}\left[\sup_{0 \leq nh
\leq T}|y^*_n|^{2p}\right]$ gives the bound of
$\mathbb{E}\left[\sup_{0 \leq nh \leq T}|y_n|^{2p}\right]$.

To bound $\mathbb{E}\left[\sup_{0 \leq t \leq
T}|\bar{y}(t)|^{2p}\right]$, we denote by $n_t$ the integer for
which $t \in [t_{n_t},t_{n_t+1})$. By definitions of (\ref{ssbe1})
and (\ref{ce1}), for $t\geq0$,
\begin{eqnarray}
\bar{y}(t)&=& y_{n_t} + (t-t_{n_t})f(y_{n_t}^*,\tilde{y}_{n_t}^*) +
g(y_{n_t}^*,\tilde{y}_{n_t}^*)\Delta w_{n_t}(t) \nonumber \\ &=&
y_{n_t} + \gamma(y_{n_t}^* - y_{n_t}) +
g(y_{n_t}^*,\tilde{y}_{n_t}^*)\Delta w_{n_t}(t) \nonumber \\ &=&
(1-\gamma)y_{n_t} + \gamma y_{n_t}^* +
g(y_{n_t}^*,\tilde{y}_{n_t}^*)\Delta w_{n_t}(t),
\end{eqnarray}
where $\gamma = (t-t_{n_t})/h <1$. Thus
\begin{eqnarray}
\mathbb{E}\left[\sup\limits_{0 \leq t \leq
T}|\bar{y}(t)|^{2p}\right] &\leq& 2^{2p-1}\left\{ \gamma
\mathbb{E}\left[\sup\limits_{0 \leq nh \leq
T}|y^*_n|^{2p}\right]+(1-\gamma) \mathbb{E}\left[\sup\limits_{0 \leq
nh \leq T}|y_n|^{2p}\right]\right. \nonumber \\ &&
\left.+\mathbb{E}\left[\sup\limits_{0 \leq t \leq
T}|g(y^*_{n_t},\tilde{y}^*_{n_t})\Delta w_{n_t}(t)|^{2p}\right]
\right\}.\label{MB1}
\end{eqnarray}
Using Doob's martingale inequality \cite[Theorem 1.3.8]{Mao97}, we
derive that
\begin{eqnarray}
\mathbb{E}\left[\sup\limits_{0 \leq t \leq
T}|g(y^*_{n_t},\tilde{y}^*_{n_t})\Delta w_{n_t}(t)|^{2p}\right] \leq
\sum_{n=0}^{N}\mathbb{E}\left[\sup\limits_{0 \leq s \leq
h}|g(y^*_{n},\tilde{y}^*_{n})\Delta w_{n}(s)|^{2p}\right] \nonumber \\
\leq \left(\frac{2p}{2p-1}\right)^{2p}\sum_{n=0}^{N}\mathbb{E}
\left[|g(y^*_{n},\tilde{y}^*_{n})\Delta
w_{n}(h)|^{2p}\right].\label{ge}
\end{eqnarray}
Thus the last term in (\ref{MB1}) is bounded by considering
(\ref{estimate3}) and bounds of $\mathbb{E}[\sup_{0 \leq nh \leq
T}|y^*_n|^{2p}]$, $\mathbb{E}[\sup_{0 \leq nh \leq T}|y_n|^{2p}]$.
Now boundedness of $\mathbb{E}[\sup_{0 \leq t \leq
T}|\bar{y}(t)|^{2p}]$ follows immediately.

\begin{lem} \label{lem6}
Under Assumption \ref{OLC}, if $(\gamma_1 + \gamma_2)h< 1$, the
implicit equation in (\ref{ssbe1}) admits a unique solution.
\end{lem}

{\it Proof.} Let $\tilde{f}(c) := f(c,\mu c + (1-\mu)b)$, then the
implicit equation (\ref{ssbe1}) takes the form as
$$
c = h \tilde{f}(c) + d = h f(c,\mu c + (1-\mu)b) + d,
$$
where at each step, $0\leq\mu<1, b, d$ are known. Observing that
\begin{eqnarray}
\langle c_1-c_2, \tilde{f}(c_1)-\tilde{f}(c_2) \rangle &=& \langle
c_1-c_2, f(c_1,\mu c_1 + (1-\mu)b)-f(c_2,\mu c_1 + (1-\mu)b) \rangle
\nonumber \\ && + \langle c_1-c_2, f(c_2,\mu c_1 +
(1-\mu)b)-f(c_2,\mu c_2 + (1-\mu)b) \rangle \nonumber \\ &\leq&
\gamma_1 |c_1-c_2|^2 + \mu \gamma_2|c_1-c_2|^2 \nonumber \\ &\leq&
(\gamma_1+\gamma_2) |c_1-c_2|^2,\nonumber
\end{eqnarray}
the assertion follows immediately from Theorem 14.2 of \cite{HW96}.
\begin{cor}
Under Assumption \ref{init},\ref{OLC}, if $(\gamma_1+\gamma_2)h<1$,
then the numerical solution produced by (\ref{ssbe1})-(\ref{ssbe2})
is well-defined and will converge to the true solution in the
mean-square sense, i.e.,
$$
\mathbb{E} \sup_{0 \leq t \leq T} \left|\bar{y}(t)-x(t)\right|^2
\rightarrow 0, \quad as \quad h \rightarrow 0.
$$
\end{cor}

{\it Proof.} Noticing that Assumption \ref{OLC} implies Assumptions
\ref{LCMC},\ref{ANMB} by Theorem \ref{EU} and Lemma \ref{MBlem}, and
taking Lemma \ref{lem6}  into consideration, the result follows
directly from Theorem \ref{ssbemain}.

\begin{rem}
We remark that the problem class satisfying condition
(\ref{initial2}) includes plenty of important models. In particular,
stochastic pantograph differential equations (see, e.g., \cite{FMW07}) with $\tau(t) = (1-q)t,
0 < q <1$ and SDDEs with constant lag fall into this class and
therefore corresponding convergence results follow immediately.
\end{rem}

\section{Mean-square stability with bounded delay}

In this section, we will investigate how SSBE shares exponential
mean-square stability of general nonlinear systems. In deterministic
case, nonlinear stability analysis of numerical methods are carried
on under a one-sided Lipschitz condition. This phenomenon has been
well studied in the deterministic case (\cite{BZ03,HW96} and
references therein) and stochastic case without delay
\cite{HMS02,HMS03,HK05,YH96,WG09a}. In what follows, we choose the
test problem satisfying  conditions (\ref{OLC1})-(\ref{OLC3}).
Moreover, we assume that variable delay is bounded , that is,  there
exists $\tau>0$, for $1\leq \kappa \in \mathbb{Z}^+, \:
0\leq\delta<1$
\begin{equation} \label{bd}
0 \leq \tau(t) \leq \tau, \quad \tau = (\kappa - \delta)h.
\end{equation}
%In \cite{WGb}, provided system (\ref{sddes1}) admits a null
%solution, exponential mean-square stability of analytical solution
%and numerical solution produced by $\theta$-Maruyama method, of
%SDDEs satisfying (\ref{bd}) is studied.
We remark that this assumption does not
impose additional restrictions on the stepsize $h$ and admits arbitrary large $h$ on
choosing $\kappa=1$ and $0\leq\delta<1$ close to 1. To begin with,
we shall first give a sufficient condition for exponential
mean-square stability of analytical solution to underlying problem.

\begin{thm} \label{ems1}
Under the conditions (\ref{OLC1}),(\ref{OLC2}),(\ref{OLC3}) and
(\ref{bd}), and with $\gamma_1, \gamma_2, \gamma_3, \gamma_4$ obeying
\begin{equation} \label{beta}
\beta := 2\gamma_1 + 2\gamma_2+\gamma_3+\gamma_4 < 0,
\end{equation}
any two solutions $x(t;\psi)$ and $y(t;\phi)$ with
$\mathbb{E}\|\psi\|^2 < \infty$ and $\mathbb{E}\|\phi\|^2 < \infty$
satisfy
$$
\mathbb{E}|x(t)-y(t)|^2 \leq \mathbb{E}\|\phi-\psi\|^2 \exp\{-\nu^+
t\},
$$
where $\nu^+ \in (0,-\beta]$ is the zero of $\mathcal {L}(\nu)=\nu +
\beta_1+\beta_2\exp\{\nu \tau\}$, with $\beta_1=2\gamma_1 +
\gamma_2+\gamma_3,\: \beta_2 =\gamma_2+\gamma_4$.
\end{thm}

{\it Proof.} By It\^{o} formula, we have
\begin{eqnarray} \label{5.12}
&&\mathbb{E}|x(t+\delta)-y(t+\delta)|^2 - \mathbb{E}|x(t)-y(t)|^2
\nonumber \\ &=& \int_t^{t+\delta}2\mathbb{E}\langle x(s)-y(s),
f(x(s),x(s-\tau(s)))-f(y(s),y(s-\tau(s)))\rangle \mbox{d}s\nonumber
\\ &&+ \int_t^{t+\delta}\mathbb{E}|g(x(s),x(s-\tau(s)))-g(y(s),y(s-\tau(s)))|^2
\mbox{d}s \allowdisplaybreaks\nonumber
\\ &\leq& (2\gamma_1+ \gamma_3)
\int_t^{t+\delta}\mathbb{E}|x(s)-y(s)|^2 \mbox{d}s +
\gamma_4\int_t^{t+\delta}\mathbb{E}|x(s-\tau(s))-y(s-\tau(s))|^2\mbox{d}s
\nonumber \\ && + 2 \int_t^{t+\delta}\mathbb{E}\langle x(s)-y(s),
f(y(s),x(s-\tau(s)))-f(y(s),y(s-\tau(s)))\rangle \mbox{d}s\nonumber \\
&\leq& \beta_1 \int_t^{t+\delta} \mathbb{E}|x(s)-y(s)|^2 \mbox{d}s+
\beta_2\int_t^{t+\delta}\sup_{r\in[s-\tau,s]}
\mathbb{E}|x(r)-y(r)|^2 \mbox{d}s.
\end{eqnarray}
Letting $u(t)=\mathbb{E}|x(t)-y(t)|^2$ and noticing that $u(t)$
exists for $t \geq -\tau$ and is continuous, we derive from
(\ref{5.12}) that
$$
\mbox{D}^+u(t) \leq \beta_1u(t)+ \beta_2\sup_{s\in[t-\tau,t]}u(s),
$$
where the upper Dini derivative $\mbox{D}^+u(t)$ is defined as
$$
\mbox{D}^+u(t):= \limsup_{\delta \rightarrow 0+}
\frac{u(t+\delta)-u(t)}{\delta}.
$$
Using Theorem 7 in \cite{BB05} leads to the desired result.
%$$
%u(t) \leq \left\{\sup_{s\in[-\tau,0]}u(s) \right \} \exp\{-\nu^+t\}.
%$$
%Substituting $u(t)$ with $\mathbb{E}|x(t)-y(t)|^2$ yields the
%desired result.

%\subsection{Mean-square stability of split-step backward Euler method}
Based on this stability result, we are going to investigate
stability of the numerical method.

\begin{thm} \label{ems3}
Under the conditions (\ref{OLC1}),(\ref{OLC2}),(\ref{OLC3}) and
(\ref{bd}), if $\beta < 0$, then for all $h > 0$, any two solutions
$X_n,Y_n$ produced by SSBE (\ref{ssbe1})-(\ref{ssbe2}) with
$\mathbb{E}\|\psi\|^2 < \infty$ and $\mathbb{E}\|\phi\|^2 < \infty$
satisfy
$$
\mathbb{E}|X_n-Y_n|^2 \leq \mathbb{E}\|\phi-\psi\|^2 \exp\{-\nu^+_h
nh\}, \quad \mbox{as} \quad n \rightarrow \infty,
$$
where $\nu^+_h >0 $ is defined as
\begin{equation}
\nu_h^+ = \frac{1}{2(\kappa+1)h}
\ln\left(\frac{1-2h\gamma_1-h\gamma_2}{1+ h \gamma_2 +h\gamma_3+
h\gamma_4}\right) > 0. \label{nuh2}
\end{equation}
\end{thm}

{\it Proof.} Under $\beta < 0$, the first part is an immediate
result from Lemma \ref{lem6}. For the second part, in order to state
conveniently, we introduce some notations
\begin{equation} \label{5.9}
W_n^* = X_n^*-Y_n^*,\: \Delta f_n^* = f(X_n^*,
\tilde{X}^*_n)-f(Y_n^*, \tilde{Y}^*_n), \: \Delta g^*_n = g(X_n^*,
\tilde{X}^*_n)-g(Y_n^*, \tilde{Y}^*_n).
\end{equation}
From (\ref{y*n2}), we have
\begin{equation} \label{5.7}
W_{n}^* = W_{n-1}^* + h \Delta f^*_{n} + \Delta g^*_{n-1} \Delta
w_{n-1}.
\end{equation}
Thus
$$
|W_n^*-h \Delta f^*_{n}|^2 = |W_{n-1}^*|^2+ 2\langle W_{n-1}^*,
\Delta g^*_{n-1} \Delta w_{n-1}\rangle + |\Delta g^*_{n-1} \Delta
w_{n-1}|^2.
$$
Taking expectation and using (\ref{OLC3}) yields
\begin{equation} \label{5.3}
\mathbb{E}|W_{n}^*|^2 - 2h \mathbb{E} \langle W_{n}^*, \Delta
f^*_{n}\rangle \leq (1+ h\gamma_3)\mathbb{E}|W_{n-1}^*|^2 +
h\gamma_4 \mathbb{E}|\tilde{X}_{n-1}^*-\tilde{Y}_{n-1}^*|^2.
\end{equation}
Now using the Cauchy-Schwarz inequality and conditions
(\ref{OLC1})-(\ref{OLC2}), we have
\begin{eqnarray}
2 \mathbb{E} \langle W_{n}^*, \Delta f^*_{n}\rangle &=& 2 \mathbb{E}
\langle W_{n}^*, f(X_{n}^*,\tilde{X}_{n}^*)
-f(Y_{n}^*,\tilde{X}_{n}^*)\rangle \nonumber
\\ && +2 \mathbb{E} \langle W_{n}^*,
f(Y_{n}^*,\tilde{X}_{n}^*)-f(Y_{n}^*,\tilde{Y}_{n}^*)\rangle \nonumber \\
&\leq& 2\gamma_1\mathbb{E}|W_{n}^*|^2 +
2\gamma_2\mathbb{E}|W_{n}^*||\tilde{X}_{n}^*-\tilde{Y}_{n}^*| \nonumber \\
&\leq& (2\gamma_1+\gamma_2)\mathbb{E}|W_{n}^*|^2 + \gamma_2
\mathbb{E}|\tilde{X}_{n}^*-\tilde{Y}_{n}^*|^2.\nonumber
\end{eqnarray}
Inserting it into (\ref{5.3}) gives
\begin{eqnarray}
(1-2h\gamma_1-h\gamma_2)\mathbb{E}|X_n^*-Y_n^*|^2
\leq(1+ h\gamma_3)\mathbb{E}|X_{n-1}^*-Y_{n-1}^*|^2 \nonumber \\
+ h\gamma_4 \mathbb{E}|\tilde{X}_{n-1}^*-\tilde{Y}_{n-1}^*|^2 +
h\gamma_2 \mathbb{E}|\tilde{X}_{n}^*-\tilde{Y}_{n}^*|^2. \label{5.4}
\end{eqnarray}
Here we have to consider which approach is chosen to treat memory
values on non-grid points, piecewise constant interpolation ($\mu
\equiv0$) or piecewise linear interpolation. In the latter case, let
us consider two possible cases:

$\bullet$ If $\tau(t_{n}) = \tilde{\mu} h, \: 0 \leq \tilde{\mu} <
1$, then
\begin{equation} \label{5.10}
\begin{split}
\mathbb{E}|\tilde{X}_{n}^*-\tilde{Y}_{n}^*|^2 &=
\mathbb{E}|\tilde{\mu}\tilde{X}_{n-1}^*
+(1-\tilde{\mu})\tilde{X}_{n}^*-\tilde{\mu}\tilde{Y}_{n-1}^*
-(1-\tilde{\mu})\tilde{Y}_{n}^*|^2 \\ &\leq
\tilde{\mu}\mathbb{E}|\tilde{X}_{n-1}^*-\tilde{Y}_{n-1}^*|^2 +
(1-\tilde{\mu})\mathbb{E}|\tilde{X}_{n}^*-\tilde{Y}_{n}^*|^2.
\end{split}
\end{equation}
Inserting (\ref{5.10}), we derive from (\ref{5.4}) that
\begin{equation}\nonumber
\begin{split}
&[1-2h\gamma_1-(2-\tilde{\mu})h\gamma_2]\mathbb{E}|X_n^*-Y_n^*|^2
\\ &\leq(1+ h\gamma_3+\tilde{\mu} h
\gamma_2)\mathbb{E}|X_{n-1}^*-Y_{n-1}^*|^2 + h\gamma_4
\mathbb{E}|\tilde{X}_{n-1}^*-\tilde{Y}_{n-1}^*|^2.
\end{split}
\end{equation}
Hence using the fact $\beta < 0$ in (\ref{beta}) gives
\begin{eqnarray}
\mathbb{E}|X_n^*-Y_n^*|^2 &\leq& \frac{1+ h\gamma_3+\tilde{\mu}
h\gamma_2+h\gamma_4}{1-2h\gamma_1-(2-\tilde{\mu})h\gamma_2}
\max_{n-\kappa-1\leq
i\leq n-1}\mathbb{E}|X_i^*-Y_i^*|^2 \allowdisplaybreaks \nonumber \\
&\leq& \frac{1+ h \gamma_2 +
h\gamma_3+h\gamma_4}{1-2h\gamma_1-h\gamma_2} \max_{n-\kappa-1\leq
i\leq n-1}\mathbb{E}|X_i^*-Y_i^*|^2. \label{5.5}
\end{eqnarray}

$\bullet$ If $\tau(t_{n}) \geq h$, it follows from (\ref{5.4}) and
$\beta <0$ that
\begin{equation}
\mathbb{E}|X_n^*-Y_n^*|^2 \leq \frac{1+ h \gamma_2 +
h\gamma_3+h\gamma_4}{1-2h\gamma_1-h\gamma_2} \max_{n-\kappa-1\leq
i\leq n-1}\mathbb{E}|X_i^*-Y_i^*|^2. \label{5.6}
\end{equation}
Therefore, it is always true that inequality (\ref{5.6}) holds for
piecewise linear interpolation case. Obviously (\ref{5.6}) also
stands in piecewise constant interpolation case.

Further, from (\ref{ssbe1}) one sees
\begin{equation} \label{5.8}
|X_0^*-Y_0^* -h(f(X_0^*, \tilde{X}^*_0)-f(Y_0^*, \tilde{Y}^*_0))|^2
= |X_0-Y_0|^2 . \nonumber
\end{equation}
Using a similar approach as before, one can derive
\begin{equation} \label{5.11}
\mathbb{E}|X_0^*-Y_0^*|^2 \leq
\frac{1+h\gamma_2}{1-2h\gamma_1-h\gamma_2}\mathbb{E}\|\psi-\phi\|^2
\leq\mathbb{E}\|\psi-\phi\|^2.
\end{equation}
Denote
\begin{equation} \label{betah}
\beta_h :=\frac{1+ h \gamma_2 +h\gamma_3+
h\gamma_4}{1-2h\gamma_1-h\gamma_2}.
\end{equation}
Noticing that $\beta < 0$, one can readily derive $0<\beta_h<1$, we can
deduce from (\ref{5.6}) and (\ref{5.11}) that
\begin{eqnarray}
  \mathbb{E}|X^*_{n-1}-Y^*_{n-1}|^2 \leq \beta_h^
  {\lfloor\frac{n-2}{\kappa+1}\rfloor+1}E\|\psi-\phi\|^2 \leq\beta_h^
  {\frac{n-2}{\kappa+1}}E\|\psi-\phi\|^2. \nonumber
\end{eqnarray}
Here $\lfloor x \rfloor$ denotes the greatest integer less than or
equal to $x$. \\
Finally from (\ref{ssbe2}), we have for large $n$ such that
$\frac{\gamma_3+\gamma_4}{n}+\frac{n-\kappa-2}{nh(\kappa+1)}\ln
\beta_h < \frac{\ln \beta_h}{2(\kappa+1)h}$
\begin{eqnarray}
\mathbb{E}|X_{n}-Y_{n}|^2 &\leq&
(1+h\gamma_3)\mathbb{E}|X^*_{n-1}-Y^*_{n-1}|^2 +
h\gamma_4\mathbb{E}|\tilde{X}^*_{n-1}-\tilde{Y}^*_{n-1}|^2 \nonumber \\
&\leq& (1+h\gamma_3)\beta_h^
  {\frac{n-2}{\kappa+1}}E\|\psi-\phi\|^2 + h\gamma_4 \beta_h^
  {\frac{n-\kappa-2}{\kappa+1}}E\|\psi-\phi\|^2 \nonumber \\
&\leq& \mbox{e}^{(\gamma_3+\gamma_4)h}\beta_h^
  {\frac{n-\kappa-2}{\kappa+1}}E\|\psi-\phi\|^2
  \nonumber \\
&\leq& \mathbb{E}\|\phi-\psi\|^2 \exp\{-\nu^+_h nh\},
\end{eqnarray}
where $\nu^+_h$ is defined as in (\ref{nuh2}).

The stability result indicates that the method
(\ref{ssbe1})-(\ref{ssbe2}) can well reproduce long-time stability
of the continuous system satisfying conditions stated in Theorem
\ref{ems1}. Note that the exponential mean-square stability under
non-global Lipschitz conditions has been studied in \cite{HMS03} in
the case of nonlinear SDEs without delay. The preceding results can
be regarded as an extension of those in \cite{HMS03} to delay case.

%Thus this work extends nonlinear stability theory in the
%deterministic case, as well as stochastic case without delay.

\section{Mean-square linear stability}
\label{linear_MS}

Although the main focus of this work is on nonlinear SDDEs, in this
section we show that the SSBE (\ref{ssbe1})-(\ref{ssbe2}) has a very
desirable linear stability property. Hence, we consider the scalar,
linear test equation \cite{Liu,ZGH09} given by
\begin{equation} \label{lineartest}
\mbox{d} x(t) = (a x(t) + b x(t-\tau)) \mbox{d}t + (c x(t) + d
x(t-\tau)) \mbox{d}w(t).
\end{equation}
Note that (\ref{lineartest}) is a special case of (\ref{sddes1})
with $\tau(t) = \tau$, and satisfies conditions
(\ref{OLC1})-(\ref{OLC3}) with
$$
\gamma_1 = a,\quad \gamma_2 = |b|, \quad\gamma_3 = c^2+|cd|,\quad
\gamma_4=d^2+|cd|.
$$
By Theorem \ref{ems1}, (\ref{lineartest}) is mean-square stable if
\begin{equation} \label{linearMS}
a < -|b| - \frac{1}{2}(|c|+|d|)^2.
\end{equation}
For constraint stepsize $h = \tau/\kappa, 1\leq \kappa \in
\mathbb{Z}^+$, i.e., $\delta = 0$ in (\ref{bd}), the SSBE proposed
in our work applied to (\ref{lineartest}) produces
\begin{eqnarray} \label{SSBEW}
\left\{
    \begin{array}{ll}  y_n^* &=  y_n + h[ay_n^* + b y_{n-\kappa}^*], \\
     y_{n+1} &=  y_n^* + [cy_n^* + d y_{n-\kappa}^*]\Delta w_n.
    \end{array} \right.
\end{eqnarray}
In \cite{ZGH09}, the authors constructed a different SSBE for the
linear test equation (\ref{lineartest}) and their method applied to
(\ref{lineartest}) reads
\begin{eqnarray} \label{SSBEZ}
\left\{
    \begin{array}{ll}  z_n^* &=  z_n + h[az_n^* + b z_{n-\kappa+1}], \\
     z_{n+1} &=  z_n^* + [cz_n^* + d z_{n-\kappa+1}]\Delta w_n.
    \end{array} \right.
\end{eqnarray}
The stability results there  \cite[Theorem 4.1]{ZGH09} indicate that
under (\ref{linearMS}) the method (\ref{SSBEZ}) can only preserve
mean-square stability of (\ref{lineartest}) with stepsize
restrictions, but the new scheme (\ref{SSBEW}) exhibits a better
stability property.
\begin{cor}
For the linear equation (\ref{lineartest}), if (\ref{linearMS})
holds, then the SSBE (\ref{SSBEW}) is mean-square stable for any
stepsize $h = \tau/\kappa, 1\leq \kappa \in \mathbb{Z}^+$.
\end{cor}

{\it Proof.} The assertion readily follows from Theorem \ref{ems3}.

Apparently, the SSBE (\ref{SSBEW}) achieves an advantage over
(\ref{SSBEZ}) in stability property that the SSBE (\ref{SSBEW}) is able to inherit
stability of (\ref{lineartest}) for any
stepsize $h = \tau/\kappa, 1\leq \kappa \in \mathbb{Z}^+$. If one drops the
stepsize restriction $h=\frac{\tau}{\kappa}, \kappa \in
\mathbb{Z}^+$ and allow for arbitrary stepsize $h>0$, one can arrive at a sharper stability result from Theorem \ref{ems3}.

\begin{cor}
For the linear equation (\ref{lineartest}), if (\ref{linearMS})
holds, then the SSBE(\ref{ssbe1})-(\ref{ssbe2}) is mean-square stable for any
stepsize $h >0$.
\end{cor}

\section{Numerical experiments}

In this section we give several numerical examples to illustrate
intuitively the strong convergence and the mean-square stability
obtained in previous sections.

\subsection{A linear example}
The first test equation is a linear It\^{o} SDDE
\begin{eqnarray} \label{LSSDE}
\left\{
    \begin{array}{ll}  \mbox{d} x(t) = (a x(t) + b x(t-1)) \mbox{d}t + (c x(t) + d x(t-1))\mbox{d}w(t), \\
     x(t) =  0.5, \quad t \in [-1,0].
    \end{array} \right.
\end{eqnarray}
Denoting $y_N^{(i)}$ as the numerical approximation to
$x^{(i)}(t_N)$ at end point $t_N$ in the $i$-th simulation of all
$M$ simulations, we approximate means of absolute errors $\epsilon$
as
$$
\epsilon=\frac{1}{M}\sum_{i=1}^{M}|y_N^{(i)}-y^{(i)}(t_N)|.
$$
In our experiments, we use the SSBE (\ref{SSBEW}) to compute an
"exact solution" with small stepsize $h=2^{-12}$ and $M=5000$. We
choose two sets of parameters as follows

$\bullet$ Example I: $a=-2, b=1, c=d=0.5;$

$\bullet$ Example II: $a=-6, b=3, c=d=1.$

$\bullet$ Example III: $a=-20, b=12, c=2, d=1.$

\begin{figure}[hpt]
         \centering
         \includegraphics[width=2.5in,height=1.8in]{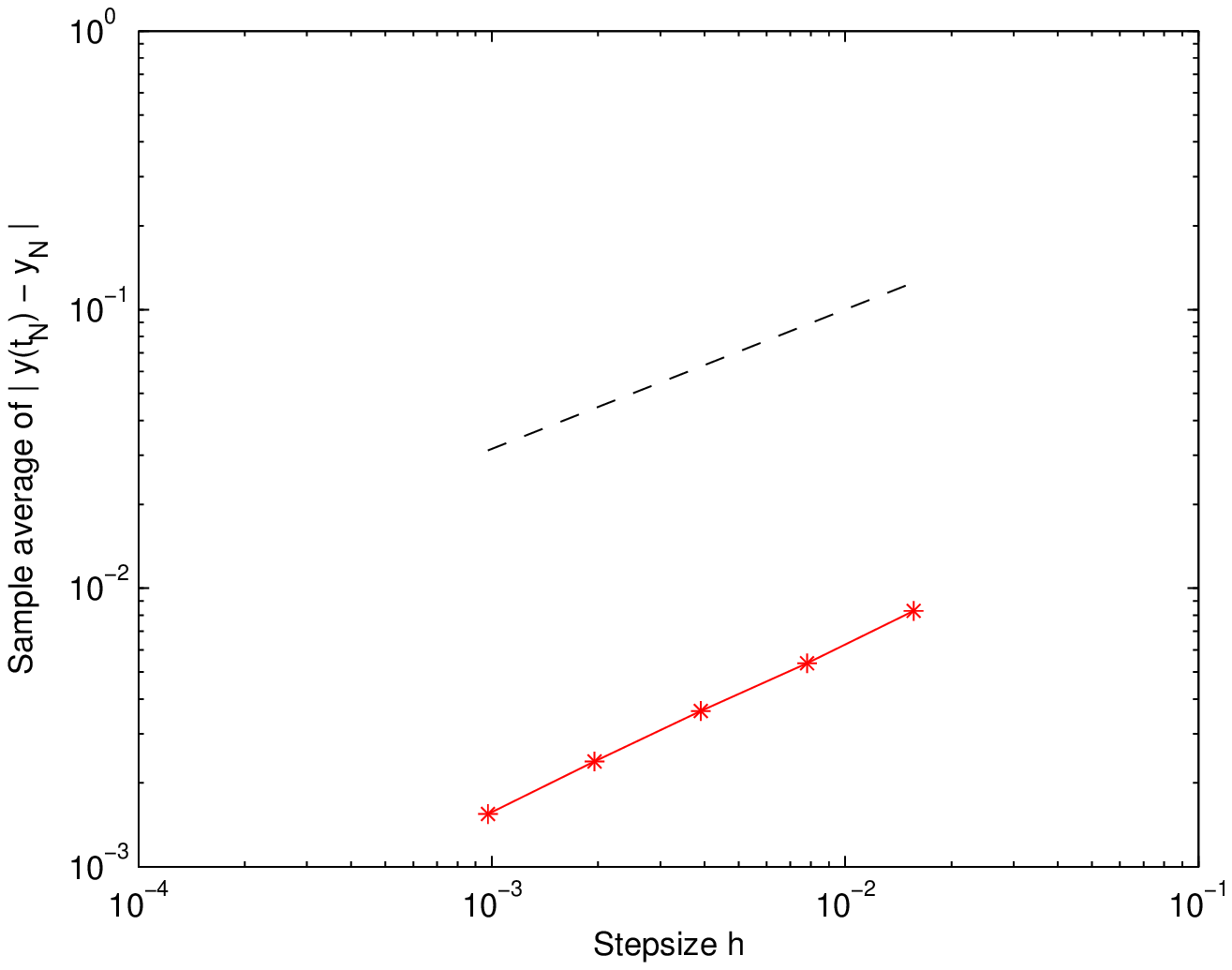}
         \includegraphics[width=2.5in,height=1.8in]{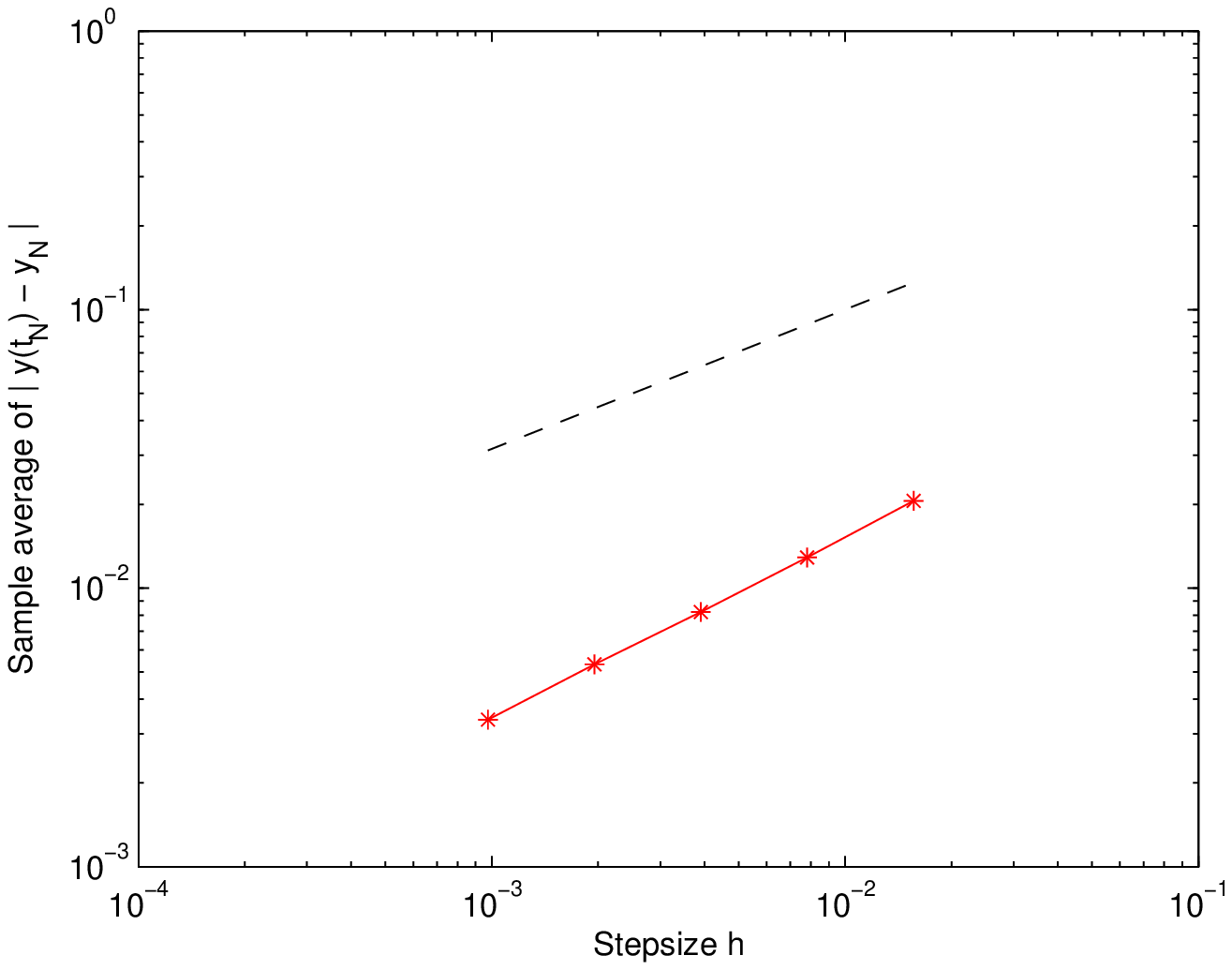}
         \caption{$\log \epsilon$ with $t_N=1$ versus $\log h$
         for Example I (left) and Example II (right).}
         \label{1}
\end{figure}

\begin{table}[h]
\begin{center} \footnotesize
\caption{Numerical results for Example II and III with $t_N=8$. }
\label{table1}
\begin{tabular*}{16cm}{@{\extracolsep{\fill}}ccccccc}
\hline \\&&\quad\quad Example II
&&& \quad\quad Example III\\
\cmidrule(r){2-4} \cmidrule(l){5-7} $h$& EM & SSBE (\ref{SSBEZ}) &
SSBE (\ref{SSBEW})
& EM & SSBE (\ref{SSBEZ}) & SSBE (\ref{SSBEW})\\
\hline\\
$2^{-7}$ &0.0008 & 0.0011& 0.0008 & 0.0014 & 0.0020& 0.0014  \\
$2^{-6}$  &0.0013 & 0.0016& 0.0013&0.0025 &0.0036& 0.0023   \\
$2^{-5}$  &0.0021 & 0.0029& 0.0019 &0.0058 &0.0070& 0.0035   \\
$2^{-4}$  &0.0034 & 0.0058& 0.0027 &0.2744 &0.0157& 0.0053   \\
$2^{-3}$  &0.0086& 0.0148& 0.0038 &6.1598e+010 &0.0628& 0.0078 \\
\hline
\end{tabular*}
\end{center}
\end{table}

In Figure \ref{1}, computational errors $\epsilon$ versus stepsize
$h$ on a log-log scale are plotted and dashed lines of slope one
half are added. One can clearly see that SSBE (\ref{SSBEW}) for
linear test equation (\ref{LSSDE}) is convergent and has strong
order of 1/2. In Table \ref{table1}, computational errors $\epsilon$
with $t_N =8$ are presented for the well-known Euler-Maruyama method
\cite{MS03}, the SSBE method (\ref{SSBEZ}) and the improved SSBE
method (\ref{SSBEW}) in this paper. There one can find that the
improved SSBE method (\ref{SSBEW}) has the best accuracy among the
three methods. In particular, for Example III with stiffness in
drift term (i.e., $a=-20$), when the moderate stepsize $h = 1/8$ was
used, the Euler-Maruyama method becomes unstable and the two SSBE
methods still remain stable, but with the improved SSBE (\ref{SSBEW})
producing better result.

\begin{figure}[hpt]
         \centering
         \includegraphics[width=2.6in,height=1.5in]{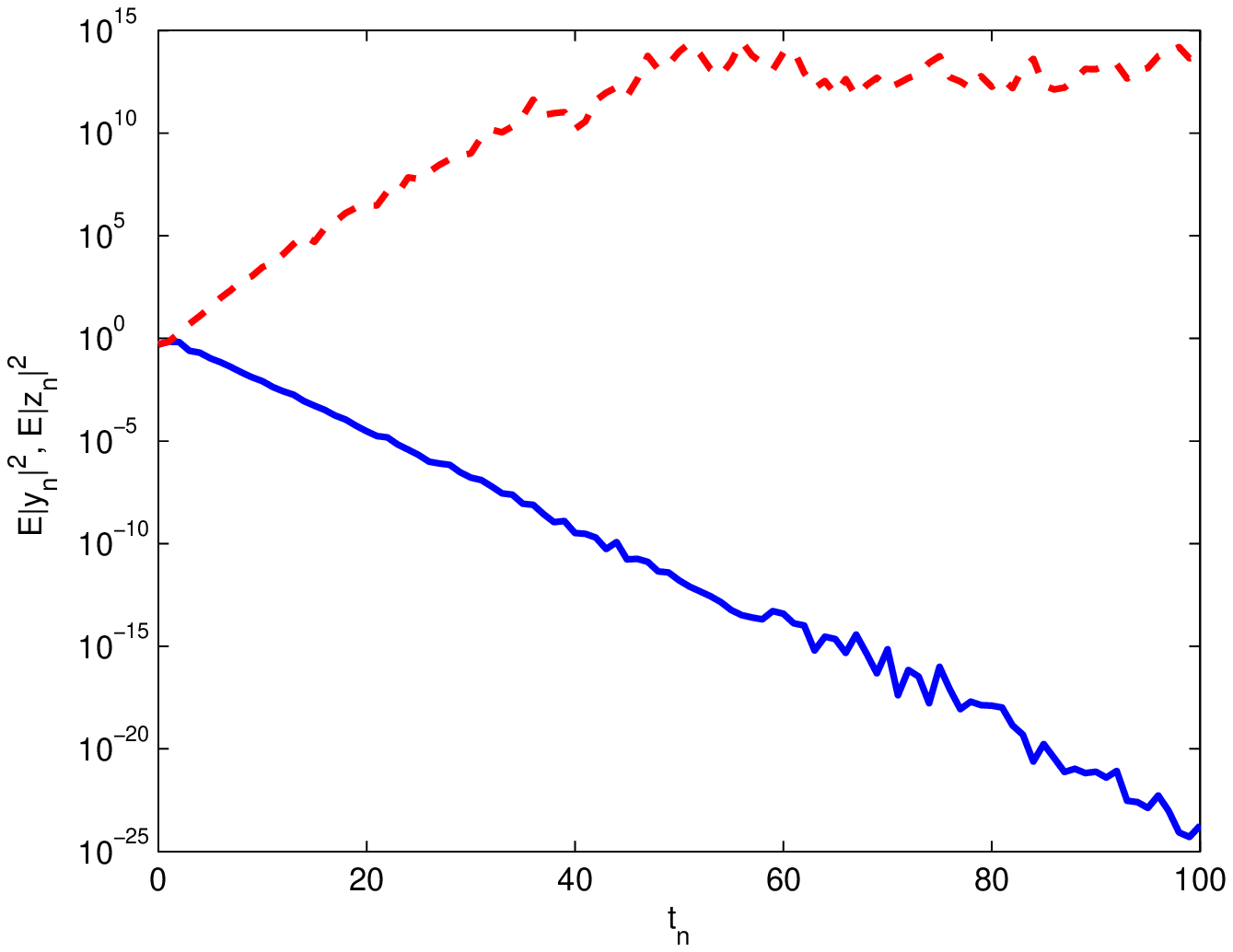}
         \includegraphics[width=2.6in,height=1.5in]{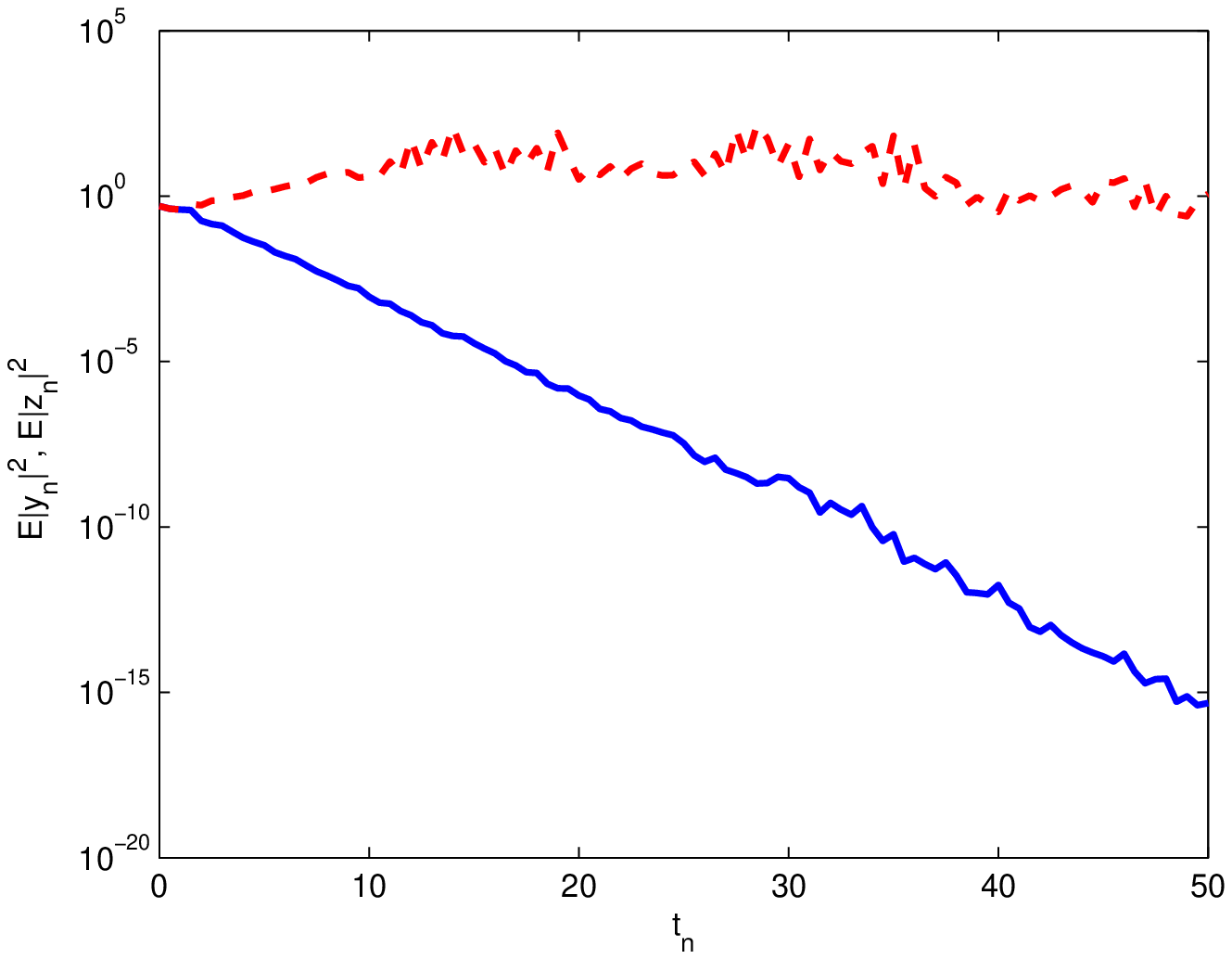}
         \includegraphics[width=2.6in,height=1.5in]{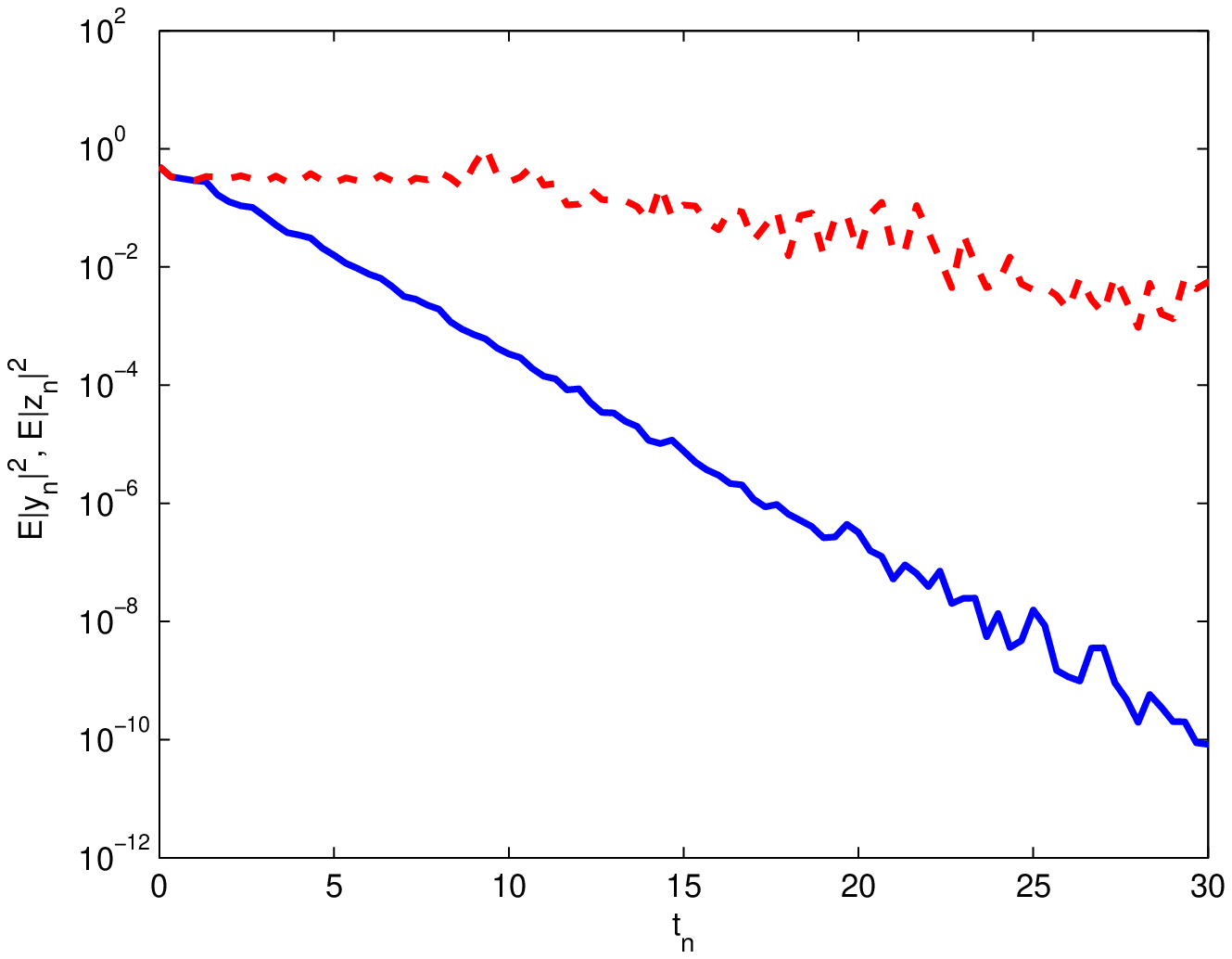}
         \includegraphics[width=2.6in,height=1.5in]{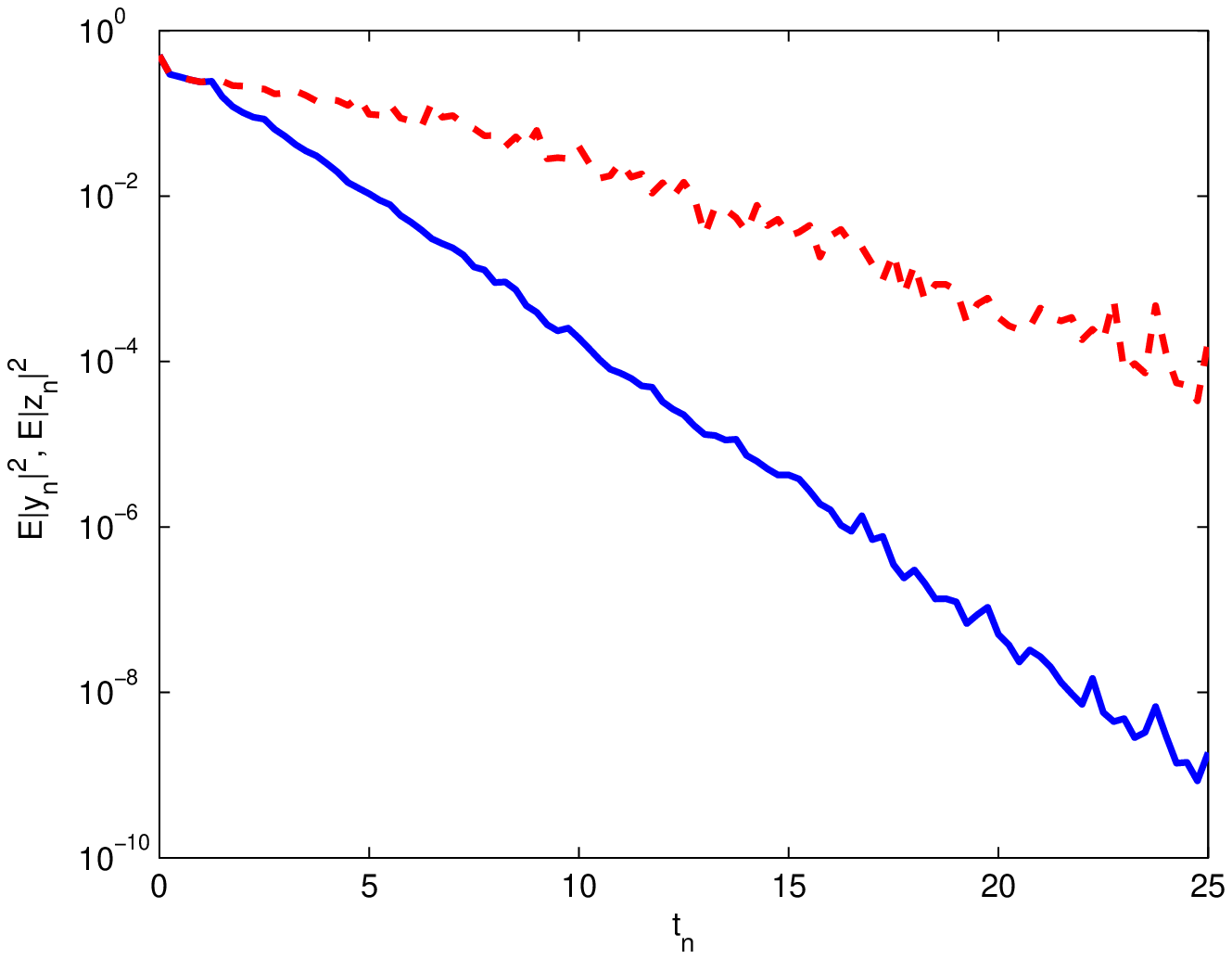}
         \caption{Simulations for (\ref{LSSDE}) with $a=-6, b=3, c=d=1$.
         Upper left: $ h= 1$, upper right: $ h= 1/2$,
         lower left: $ h= 1/3$, lower right: $ h= 1/4$.}
         \label{2}
\end{figure}

\begin{figure}[hpt]
         \centering
         \includegraphics[width=2.6in,height=1.5in]{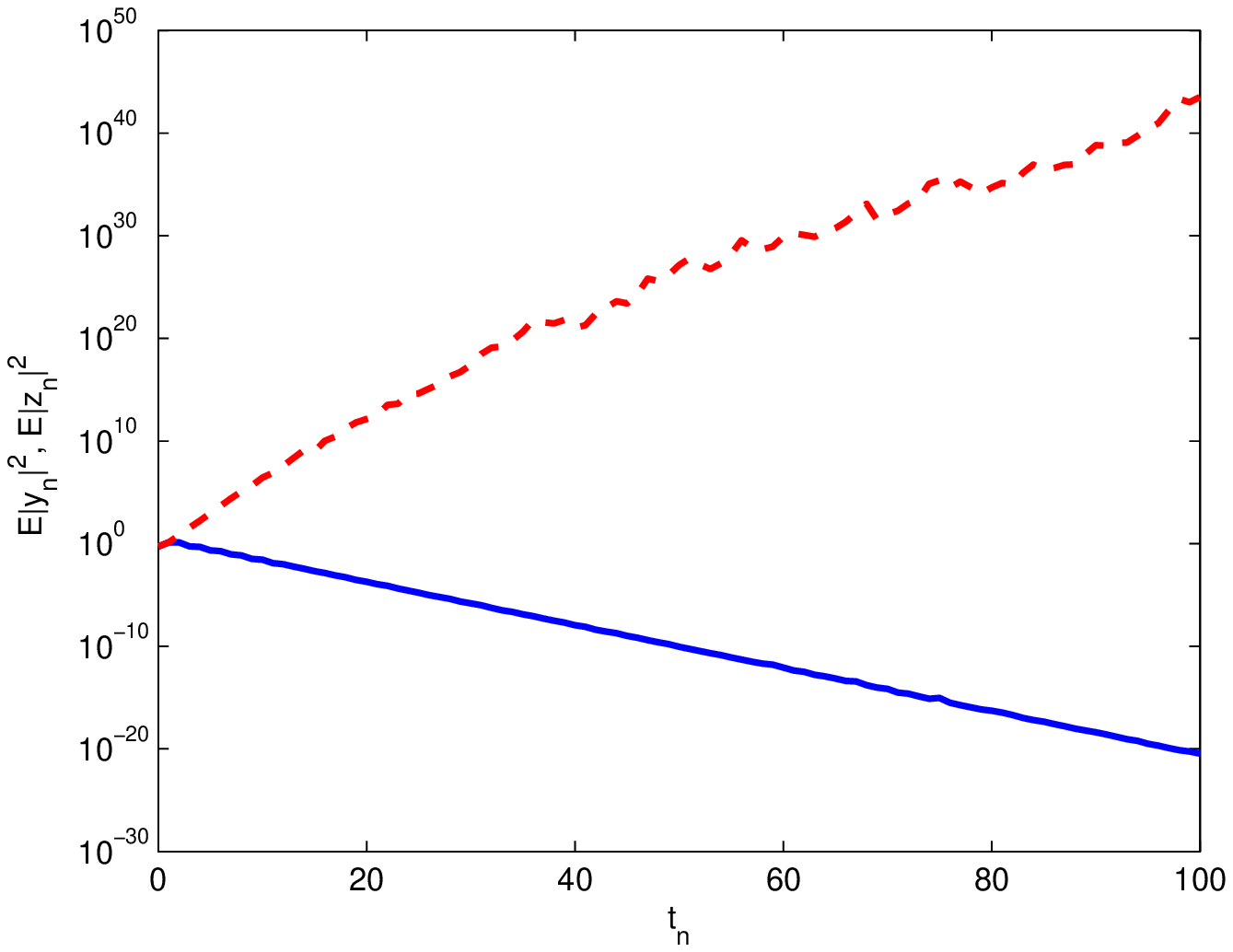}
         \includegraphics[width=2.6in,height=1.5in]{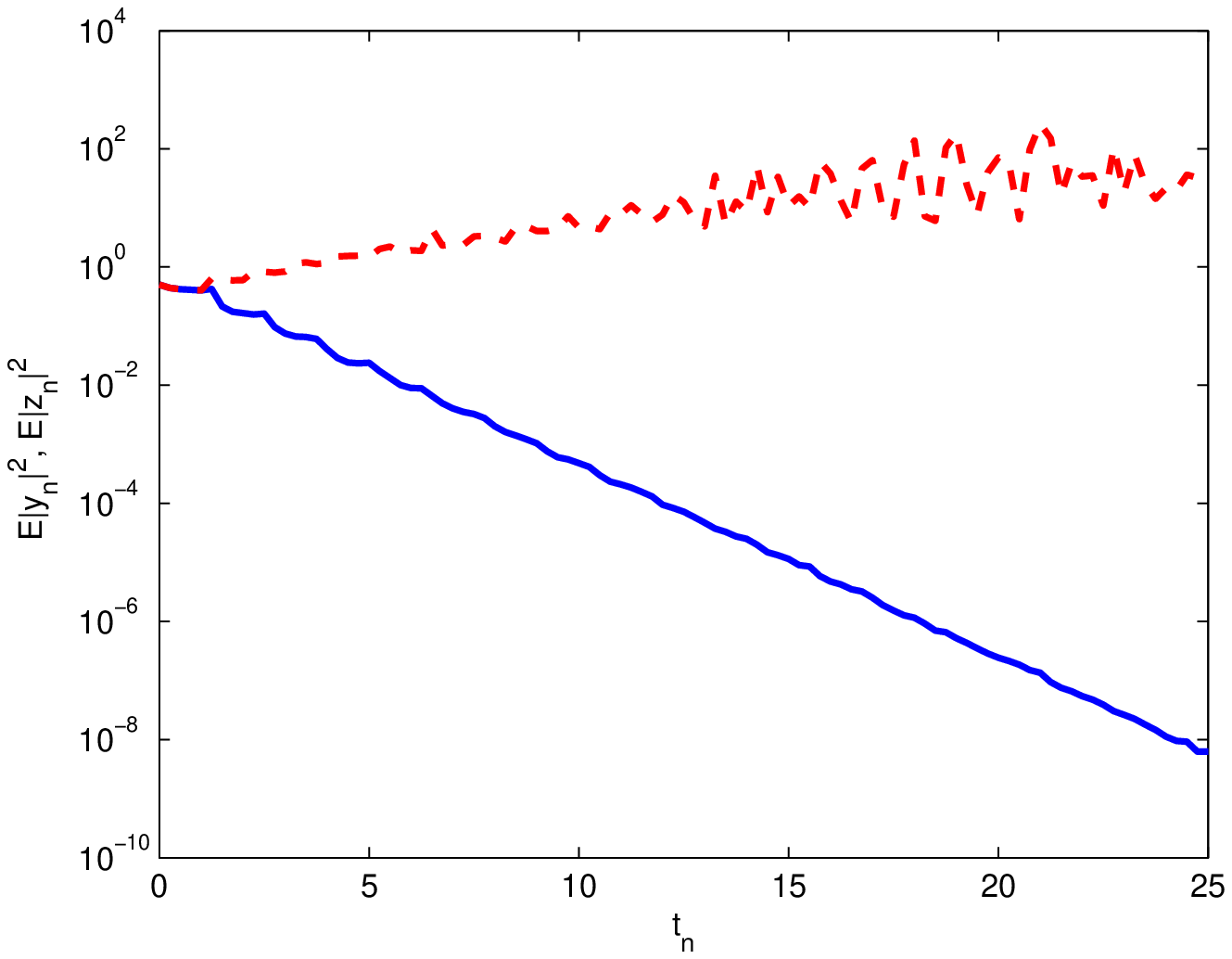}
         \includegraphics[width=2.6in,height=1.5in]{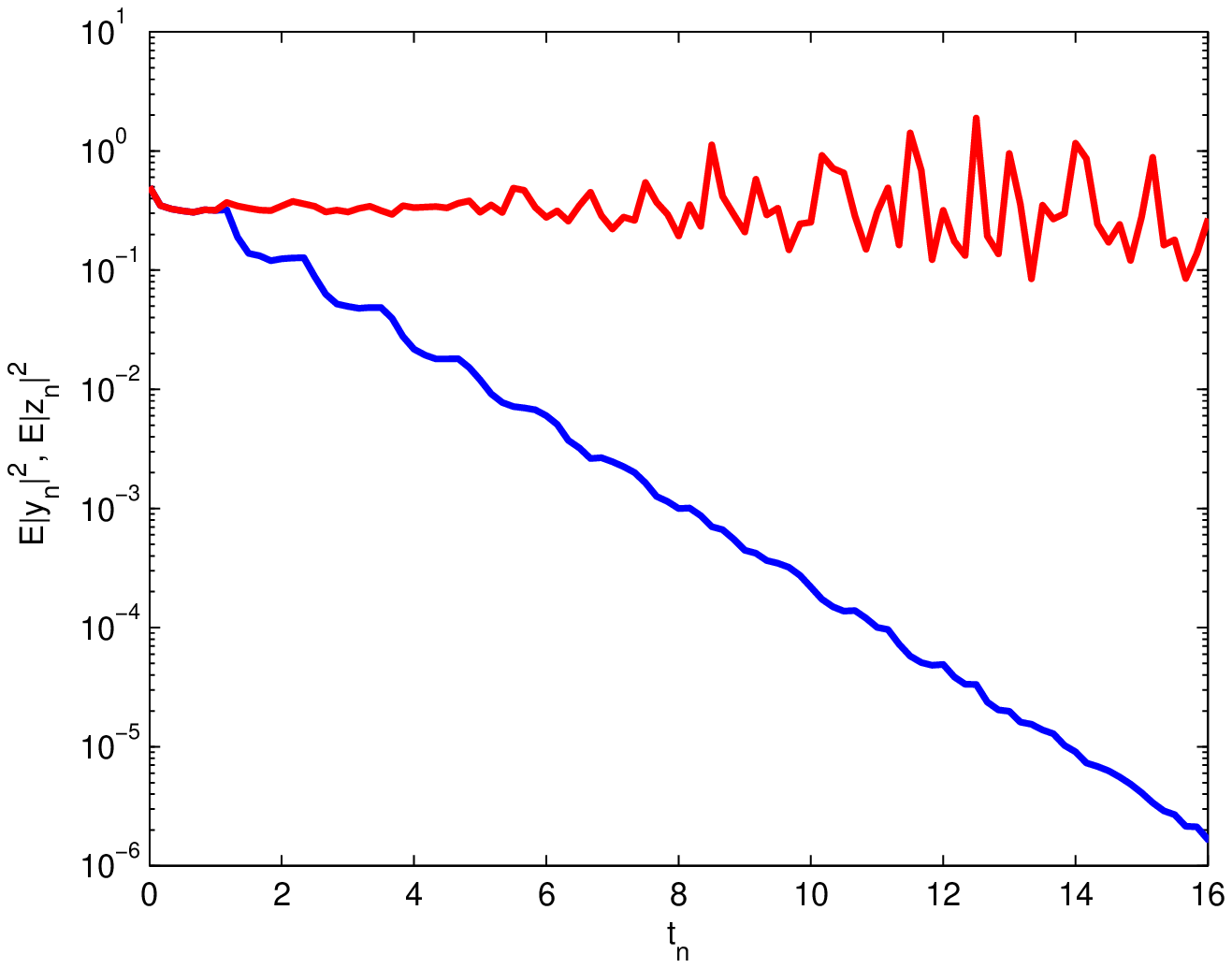}
         \includegraphics[width=2.6in,height=1.5in]{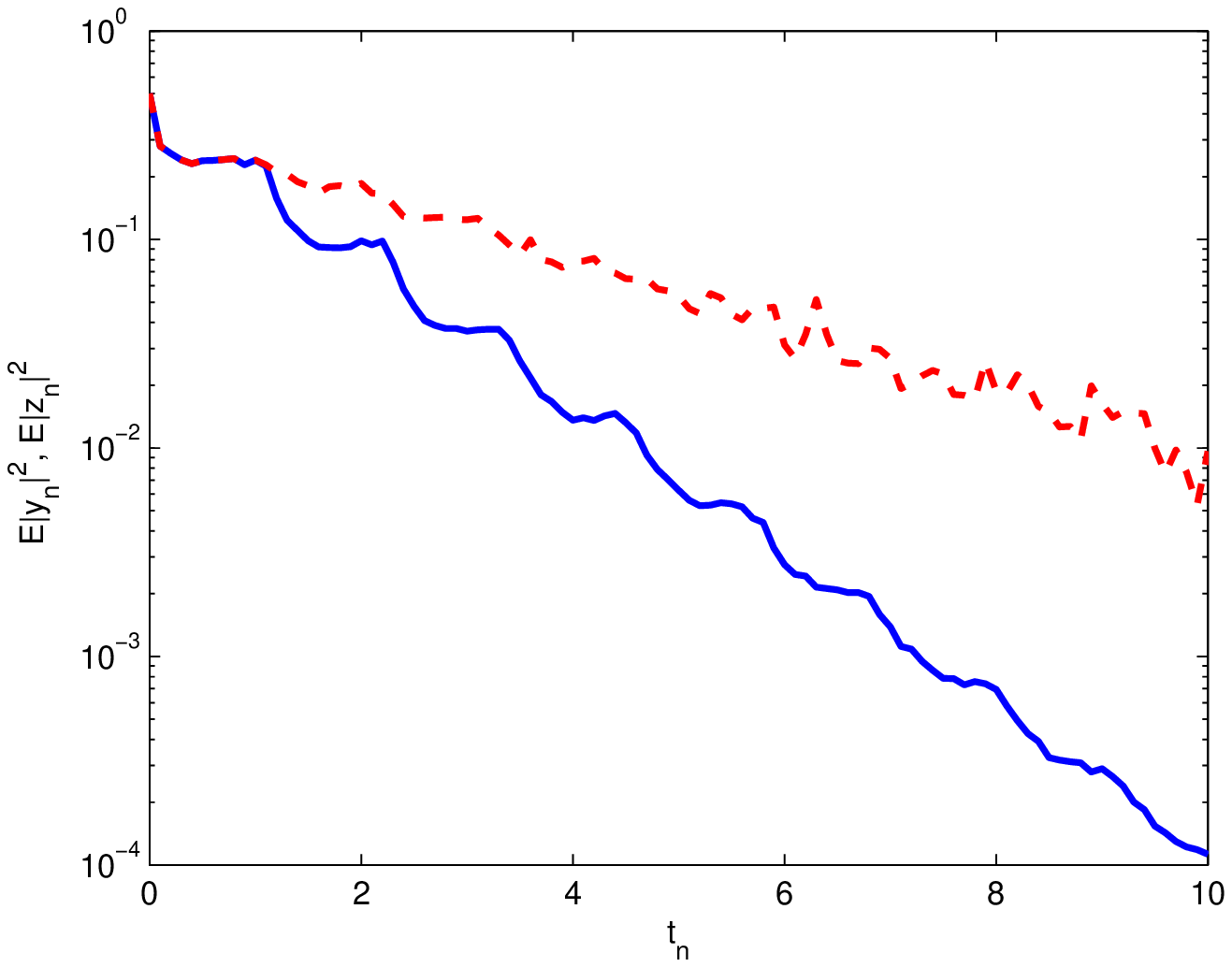}
         \caption{Simulations for (\ref{LSSDE}) with $a=-20, b=12, c=2,d=1$.
         Upper left: $ h= 1$, upper right: $ h= 1/4$,
         lower left: $ h= 1/6$, lower right: $ h= 1/10$.}
         \label{3}
\end{figure}

To compare stability property of the improved SSBE and SSBE in
\cite{ZGH09}, simulations by SSBE (\ref{SSBEW}) and (\ref{SSBEZ})
are both depicted in Figure \ref{2}, \ref{3}. There solutions
produced by (\ref{SSBEW}) and (\ref{SSBEZ}) are plotted in solid line and dashed
line, respectively. As is shown in the figures,
methods (\ref{SSBEW}) and (\ref{SSBEZ}) exhibit different stability
behavior. One can observe from Figure \ref{2} that (\ref{SSBEW}) for
Example II is mean-square stable for $h=1,1/2,1/3,1/4$. But
(\ref{SSBEZ}) is unstable for $h=1,1/2$. For Example III, the
improved SSBE (\ref{SSBEW}) is always stable for $h=1,1/4,1/6,1/10$,
but (\ref{SSBEZ}) becomes stable when the stepsize $h$ decreases to
$h=1/10$. The numerical results demonstrate that the scheme
(\ref{SSBEW}) has a greater advantage in mean-square stability than
(\ref{SSBEZ}).

%\begin{figure}[h]
%         \centering
%         \includegraphics[width=2.5in,height=1.5in]{E1.1.eps}
%         \includegraphics[width=2.5in,height=1.5in]{E1.2.eps}
%         \caption{Simulations for (\ref{LSSDE}) with $a=-2,
%         b=1, c=d=0.5$. Left : $ h= 1$, right: $ h= 1/2$. }
%         \label{2}
%\end{figure}

\subsection{A nonlinear example}

\begin{figure}[hpt]
         \centering
         \includegraphics[width=3in,height=2.5in]{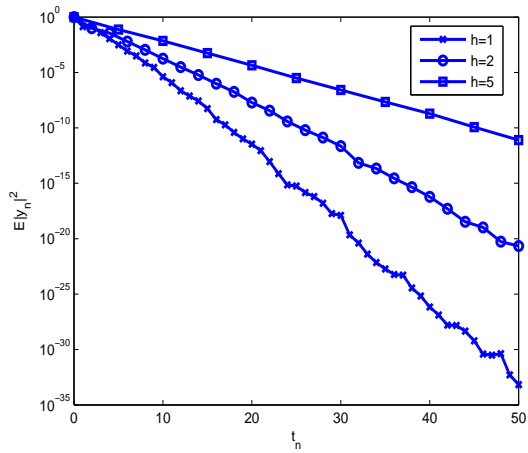}
         \caption{ Simulations for (\ref{NLSDDE}) by SSBE
         (\ref{SSBEW}) using various stepsizes.}
         \label{4}
\end{figure}

Consider a nonlinear SDDE with a time-varying delay as follows
\begin{eqnarray} \label{NLSDDE}
\left\{
    \begin{array}{ll}  d x(t) = \left[-4 x(t) - 3x^3(t)  + x(t-\tau(t))\right] dt + \left[
x(t)+x(t-\tau(t))\right] dw(t), t >0, \\
     x(t) =  1, \quad t \in [-1,0],
    \end{array} \right.
\end{eqnarray}
%\begin{equation}
%d x(t) = \left[-4 x(t) - 3x^3(t)  + x(t-\tau(t))\right] dt + \left[
%2+ x(t)+x(t-\tau(t))\right] dw(t), t >0, \label{NLSDDE}
%\end{equation}
where $\tau(t) = \frac{1}{1+t^2}$. Obviously, equation
(\ref{NLSDDE}) satisfies conditions (\ref{OLC1})-(\ref{OLC3}) in
Assumption \ref{OLC}, with $\gamma_1=-4, \gamma_2= 1, \gamma_3 =
\gamma_4 = 2$. Thus $2\gamma_1 +2 \gamma_2 + \gamma_3 + \gamma_4 =
-2 <0$ and the problem is exponentially mean-square stable. As is
shown in Figure \ref{4}, the SSBE (\ref{SSBEW}) can well reproduce
stability for quite large stepsize $h =1,2,5$. This is consistent
with our result established in Theorem \ref{ems3}.

%\section{Conclusion}

\section*{Appendix}

{\it Proof of Theorem \ref{EU}.} Since both $f$ and $g$ are locally
Lipschitz continuous, Theorem 3.2.2 of \cite{Mao94} shows that there
is a unique maximal local solution $x(t)$ on $t\in
[[0,\rho_{\infty}[[$, where the stopping time $\rho_R =
\inf\{t\geq0: |x(t)|\geq R\}$.
%Since the proof is very
%tedious, we omit it here and readers are referred to \cite{Mao97}
%for more detail.
By It\^{o}'s formula we obtain that for $t \geq 0$
\begin{eqnarray}
&&|x(t\wedge\rho_R)|^2 = |\psi(0)|^2 + 2 \int_0^{t\wedge\rho_R}
x(s)^{T}f(x(s),x(s-\tau(s))) \mbox{d}s \nonumber \\  && + 2
\int_0^{t\wedge\rho_R} x(s)^{T}g(x(s),x(s-\tau(s))) \mbox{d}w_s
+ \int_0^{t\wedge\rho_R} |g(x(s),x(s-\tau(s)))|^2 \mbox{d}s \allowdisplaybreaks \nonumber \\
&&\leq |\psi(0)|^2 + 3 K\int_0^{t\wedge\rho_R}
(1+|x(s)|^2+|x(s-\tau(s))|^2 )\mbox{d}s \nonumber \\  && + 2
\int_0^{t\wedge\rho_R} x(s)^{T}g(x(s),x(s-\tau(s))) \mbox{d}w_s,
\label{Ito1}
\end{eqnarray}
where the condition (\ref{MC}) was used. Thus
\begin{eqnarray}
\sup_{0 \leq s \leq t}|x(s\wedge\rho_R)|^2 &\leq& |\psi(0)|^2 + 3
K\int_0^{t} (1+2\sup_{0 \leq r \leq s}|x(r\wedge\rho_R)|^2+
\|\psi\|^2)\mbox{d}s \nonumber \\  &&+ 2 \sup_{0 \leq s \leq
t}\int_0^{s\wedge\rho_R} x(r)^{T}g(x(r),x(r-\tau(r))) \mbox{d}w_r.
\label{Ito2}
\end{eqnarray}
Now, raising both sides of (\ref{Ito2}) to the power $p/2$ and using
H\"{o}lder's inequality yield
\begin{eqnarray}
&&\sup_{0 \leq s \leq t}|x(s\wedge\rho_R)|^p \leq
3^{p/2-1}\left\{|\psi(0)|^p \right.\nonumber
\\ &&\left.+ (3K)^{p/2}(3T)^{p/2-1}\int_0^t (1+2^{p/2}\sup_{0 \leq r \leq
s}|x(r\wedge\rho_R)|^p+\|\psi\|^p )\mbox{d}s\right. \nonumber
\\  &&\left.+ 2^{p/2} \sup_{0 \leq s \leq
t}\left|\int_0^{s\wedge\rho_R} x(r)^{T}g(x(r),x(r-\tau(r)))
\mbox{d}w_r\right|^{p/2}\right\}. \label{Ito4}
\end{eqnarray}
By the Burkholder-Davis-Gundy inequality \cite{Mao97}, one computes
that, with $c_1 = c_1(p,T)$,
\begin{eqnarray}
\mathbb{E}\left[\sup\limits_{0 \leq s \leq
t}|x(s\wedge\rho_R)|^{p}\right] &&\leq
c_1\left\{1+\mathbb{E}\|\psi\|^p + \int_0^t \mathbb{E}\sup_{0
\leq r \leq s}|x(r\wedge\rho_R)|^p\mbox{d}s\right. \nonumber \\
&&\left. + \mathbb{E}\left[\int_0^{t\wedge\rho_R}
|x(s)|^2|g(x(s),x(s-\tau(s)))|^2 \mbox{d}s \right]^{p/4}\right\}.
\label{Ito5}
\end{eqnarray}
Next, by an elementary inequality,
\begin{eqnarray}
&& \mathbb{E}\left[\int_0^{t\wedge\rho_R}
|x(s)|^2|g(x(s),x(s-\tau(s)))|^2 \mbox{d}s \right]^{p/4} \nonumber
\\ &\leq& \mathbb{E}\left[\sup\limits_{0 \leq s \leq
t}|x(s\wedge\rho_R)|^2\int_0^{t\wedge\rho_R}
|g(x(s),x(s-\tau(s)))|^2 \mbox{d}s \right]^{p/4} \nonumber
\\ &\leq& \frac{1}{2c_1}\mathbb{E}\left[\sup\limits_{0 \leq s \leq
t}|x(s\wedge\rho_R)|^{p}\right] +
\frac{c_1}{2}T^{p/2-1}\mathbb{E}\int_0^{t\wedge\rho_R}|g(x(s),x(s-\tau(s)))|^p
\mbox{d}s \nonumber
\\ &\leq& \frac{1}{2c_1}\mathbb{E}\left[\sup\limits_{0 \leq s \leq
t}|x(s\wedge\rho_R)|^{p}\right]+\frac{c_1}{2}(3T)^{p/2-1}K^{p/2}\int_0^t
(1+ \mathbb{E}\sup_{0 \leq r \leq
s}|x(r\wedge\rho_R)|^p+\mathbb{E}\|\psi\|^p) \mbox{d}s. \nonumber
\end{eqnarray}
Inserting it into (\ref{Ito5}), for proper constants $c_2, c_3$ we
have  that
$$
\mathbb{E}\sup_{0 \leq s \leq t}|x(s\wedge\rho_R)|^p \leq
c_2(1+\mathbb{E}\|\psi\|^p) + c_3\int_0^{t}\mathbb{E}\sup_{0 \leq r
\leq s}|x(r\wedge\rho_R)|^p\mbox{d}s.
$$
The Gronwall inequality gives
\begin{equation}\label{Ito6}
\mathbb{E}\sup_{0 \leq s \leq T}|x(s\wedge\rho_R)|^p \leq
c_2(1+\mathbb{E}\|\psi\|^p) \mbox{e}^{c_3T}.
\end{equation}
This implies
$$
R^p \mathbb{P}\{\rho_R \leq T\}\leq c_2(1+\mathbb{E}\|\psi\|^p)
\mbox{e}^{c_3T}.
$$
Letting $R\rightarrow\infty$ leads to
$$
\lim_{R\rightarrow\infty}\mathbb{P}\{\rho_R \leq T\}=0.
$$
Since $T>0$ is arbitrary, we must have $\rho_R\rightarrow\infty$
a.s. and hence $\rho_{\infty}=\infty$ a.s. The existence and
uniqueness of the global solution is justified. Finally, the desired
moment bound follows from (\ref{Ito6}) by letting
$R\rightarrow\infty$ and setting $C=c_2e^{c_3T}$.

%\begin{acknowledgements}
%If you'd like to thank anyone, place your comments here
%and remove the percent signs.
%\end{acknowledgements}

%\section*{Appendix. Proof of Theorem \ref{EU}.}


\begin{thebibliography}{10}

\bibitem{BB00}  C.T.H. Baker, E. Buckwar, Numerical analysis of explicit one-step
methods  for stochastic delay differential equations, LMS J. Comput.
Math., 3(2000), pp.315-335.

\bibitem{BB05} C.T.H.  Baker,  E.  Buckwar,  Exponential  stability  in
pth mean  of  solutions, and of convergent  Euler-type solutions, to
stochastic  delay differential equations, J. Comput.  Appl. Math.
,184 (2) (2005), pp.404-427.

\bibitem{BZ03} A.Bellen, M.Zennaro, Numerical Methods for Delay Differential Equations,
Oxford University Press, Oxford, 2003.

\bibitem{BBT04} K.Burrage, P.M.Burrage, T.Tian, Numerical methods for strong
solutions of stochastic differential equations: an overview,
Proceedings: Mathematical, Physical and Engineering, Royal Society
of London 460(2004), pp.373-402.

\bibitem{FMW07} Z.Fan, M.Liu, W.Cao, Existence and uniqueness of the solutions and
convergence of semi-implicit Euler methods for stochastic pantograph
equations, J. Math. Anal. Appl., 325 (2007), pp.1142-1159.

\bibitem{HW96} E.Hairer, G.Wanner, Solving Ordinary Differential Equations II: Stiff and Differential-
Algebraic Problems, Springer-Verlag, Berlin, second ed., 1996.

\bibitem{HMS02} D.J.Higham, X.Mao, A.M.Stuart, Strong convergence of
Euler-type methods for non-linear stochastic diffrential equations,
SIAM J. Numer. Anal. 40(2002), pp.1041-1063.

\bibitem{HMS03} D.J.Higham, X.Mao, A.M.Stuart, Exponential mean-square stability
of numerical solutions to stochastic differential equations, LMS J.
Comput. Math., 6 (2003), pp.297-313.

\bibitem{HK05} D.J.Higham, P.E.Kloeden, Numerical methods for nonlinear
stochastic differential equations with jumps, Numer. Math.
101(2005), pp.101-119.

\bibitem{HK06} D.J.Higham, P.E.Kloeden, Convergence and stability of implicit
methods for jump-diffusion systems, Int. J. Numer. Anal. Model., 3
(2006) 125-140.

\bibitem{HK07} D.J.Higham, P.E.Kloeden, Strong convergence rates for backward Euler on a
class of nonlinear jump-diffusion problems, J. Comput. Appl. Math., 205 (2007) 949-956.

\bibitem{YH96} Y.Hu, Semi-implicit Euler-Maruyama scheme for stiff stochastic
equations, in Stochastic Analysis and Related Topics V: The Silvri
Workshop, Progr. Probab. 38, H. Koerezlioglu, ed., Birkhauser,
Boston, 1996, pp.183-202.

\bibitem{JKN09} A.Jentzen, P.E.Kloeden, A.Neuenkirch,  Pathwise approximation of stochastic differential equations on domains: higher order convergence rates without global Lipschitz coefficients, Numer. Math. 112, 1 (2009), 41-64.

\bibitem{KP92} P.E.Kloeden, E.Platen, Numerical Solution of Stochastic
Differential Equations, Springer, Berlin, 1992.

\bibitem{Liu} M.Liu, W.Cao, Z.Fan, Convergence and stability of the
semi-implicit Euler method for a linear stochastic differential
delay equation, J. Comput. Appl. Math., 170 (2004), pp.255-268.

\bibitem{Mao94} X.Mao, Exponential Stability of Stochastic Differential
Equations, Marcel Dekker, NewYork, 1994.

\bibitem{Mao97} X.Mao, Stochastic Differential Equations and Applications, Horwood,
New York, 1997.

\bibitem{MS03} X.Mao, S. Sabanis, Numerical solutions of stochastic differential
delay equations under local Lipschitz condition, J. Comput. Appl.
Math., 151 (2003), pp.215-227.

\bibitem{Mao07} X.Mao, Exponential stability of equidistant Euler-Maruyama
approximations of stochastic differential delay equations, J.Comput.
Appl. Math., 200 (2007), pp.297 - 316.

\bibitem{GM05} G.N.Milstein, M.V. Tretyakov, Numerical integration of
stochastic diffrential equations with nonglobally Lipschitz
coeffients, SIAM J. Numer. Anal., 43 (3)(2005), pp.1139-1154.

\bibitem{SM96} Y.Saito, T.Mitsui, Stability analysis of numerical
schemes for stochastic differential equations, SIAM J. Numer. Anal.,
33 (1996), pp.2254-2267.

\bibitem{TL89} L.Torelli, Stability of numerical methods for delay differential
equations, J. Comput. Appl. Math. 25(1989), pp.15-26.

\bibitem{ZGH09} H.Zhang, S.Gan and L.Hu, The split-step backward Euler method
for linear stochastic delay differential equations, J. Comput. Appl.
Math. 225 (2009), pp.558-568.

\bibitem{WG09a} X.Wang, S.Gan,  Compensated stochastic theta methods for stochastic
differential equations with jumps, Appl. Numer. Math. 60 (2010),
pp.877-887.

\end{thebibliography}
\end{document}